\documentclass[12pt]{article}

\usepackage{amsmath,amsbsy,amsfonts,amssymb,
amsthm,dsfont,fullpage,units,bm,mathrsfs,bbm,
empheq,graphicx,subcaption,algorithm,algorithmic,
mathtools,color,cases,tikz,multirow,hyperref,
comment,enumerate,authblk}

\hypersetup{
  colorlinks = true,
  citecolor = {blue},
  citebordercolor = {white}
}
\usepackage[top=1.in, bottom=1in, left=1.1in, right=1.1in]{geometry}

\newcommand{\bi}{\begin{itemize}}
\newcommand{\ei}{\end{itemize}}

\newtheorem{theorem}{Theorem}

\newtheorem{lemma}{Lemma}

\theoremstyle{definition}

\DeclareMathOperator{\Tr}{Tr}

\newcommand{\R}{\mathbb{R}}
\newcommand{\Z}{\mathbb{Z}}

\renewcommand{\P}{\mathbb{P}}
\newcommand{\EE}{\mathbb{E}}
\newcommand{\E}{\mathbb{E}}

\newcommand{\N}{\mathbb{N}}

\newcommand{\cH}{\mathcal{H}}

\newcommand{\bx}{{x}}
\newcommand{\bX}{{X}}

\newcommand{\bz}{\bm{z}}

\newcommand{\eps}{\varepsilon}
\newcommand{\be}{\begin{equation}}
\newcommand{\ee}{\end{equation}}

\newcommand{\transpose}{^{\operatorname{T}}}
\newcommand{\rmd}{\mathrm{d}}

\newcommand{\del}{\partial}

\begin{document}

\title{Algorithms for Solving High Dimensional PDEs: 
From Nonlinear Monte Carlo to Machine Learning}

\author[1,2]{Weinan E}
\author[1]{Jiequn Han}
\author[3]{Arnulf Jentzen} 
\affil[1]{\normalsize Department of Mathematics, Princeton University}
\affil[2]{\normalsize Program in Applied and Computational Mathematics, Princeton University}
\affil[3]{\normalsize Faculty of Mathematics and Computer Science, University of
M\"{u}nster}

\date{\normalsize \today}

\maketitle
\begin{abstract}
In recent years, tremendous progress has been made on numerical algorithms for solving partial differential equations (PDEs) in a very high dimension, using ideas from either nonlinear (multilevel) Monte Carlo or deep learning. They are potentially free of the curse of dimensionality for many different applications and have been proven to be so in the case of some nonlinear Monte Carlo methods for nonlinear parabolic PDEs.

In this paper, we review these numerical and theoretical advances. In addition to algorithms based on  stochastic reformulations of the original problem, such as the multilevel  Picard iteration and the Deep BSDE method,  we also discuss algorithms based on the more traditional Ritz, Galerkin, and least square formulations. We hope to demonstrate to the reader that studying PDEs as well as control and variational problems in very high dimensions might very well be among the most promising new directions in mathematics and scientific computing in the near future.
\end{abstract}

{
  \hypersetup{linkcolor=black}
  \tableofcontents
}

\section{Introduction}
\label{sec:intro}

The mathematical models for many problems around us are in the form of partial differential equations (PDEs) in high dimensions.
Notable examples include:
\bi
  \item The Hamilton-Jacobi-Bellman (HJB) equation in  control theory
  \begin{equation}
    \frac{ \partial u }{ \partial t } + H( \bx, \nabla_x u ) = 0 .
  \end{equation}
  Here  the dimensionality is the dimensionality of the state space of the original control problem. If the original control
  problem is described by a PDE, then the corresponding HJB equation is formulated in infinite dimensional space.
  \item The Black-Scholes equations for pricing financial derivatives
  \begin{equation}
    \frac{ \partial u }{ \partial t } + \frac12\sigma^2 \sum_{i=1}^d x_i^2\frac{\partial^2 u}{\partial x_i^2} + r \left< \nabla_x u , \bx \right> - r u + f(u) = 0 .
  \end{equation}
  Here the dimensionality is the number of underlying financial assets. 
  Nonlinear terms $f(u)$ may result when default risks, transaction costs, or other factors are taken into account.
  \item Many electron Schr\"odinger equation in quantum mechanics
\begin{equation}
i \frac{ \partial u }{ \partial t } =   \Delta_x u + V(\bx) u .
\end{equation}  
Here the dimensionality is roughly three times the number of electrons in the considered quantum mechanical system.
\ei

Solving these PDEs has been a  notoriously difficult problem in scientific computing and computational science,
due to the well-known {\bf curse of dimensionality (CoD)}: The computational complexity grows exponentially as a function
of the dimensionality of the problem \cite{Bellman1957}. In fact, for this reason, traditional numerical algorithms such as
finite difference and finite element methods have been limited to dealing with problems in a rather low dimension.
The use of sparse grids can extend the applicability to, say, around 10 dimensions.  But beyond that, there seemed to be
little hope except for special PDEs.

We are interested in PDEs and algorithms that do not suffer from CoD, i.e., 
algorithms whose computational complexity scales algebraically with the problem dimension.
More precisely, 
to reach some error tolerance $ \eps > 0 $, the computational cost should be no more than
\begin{equation}
c(d, \eps) \sim C d^{\alpha} \eps^{-\beta}
\label{cost}
\end{equation}
where $ C, \alpha, \beta \geq 0 $ are absolute, dimension-independent constants. 
In particular, they do not depend on the dimension
$ d \in \N = \{ 1, 2, 3, \dots \} $. We are also interested in PDEs for which such algorithms do exist.  In fact, we are interested in developing a theory of high dimensional
PDEs based on the complexity with which the solutions can be approximated using particular schemes, such as neural network models.
We believe that such a theory should be part of the foundation for a theoretical understanding of high dimensional control theory,
reinforcement learning, and a host of other problems that will become important research topics in mathematics.
 
The golden standard for high dimensional problems is the approximation of high dimensional integrals.
Let $ g \colon [0,1]^d \to \R $ be a Lebesgue square integrable function defined on the set $ X = [0, 1]^d $ and let
\begin{equation}
\label{eq:integral_intro}
  I(g) = \int_{X} g(\bx) \, dx
  .
\end{equation}
Typical grid-based quadrature rules, such as the Trapezoidal rule and the Simpson's rule,
all suffer from  CoD. The one algorithm that does not suffer from CoD
is the Monte Carlo algorithm which works as follows.
Let $ \bx_j $, $ j \in \{ 1, 2, \dots, n \} $, be 
independent, continuous uniformly distributed random samples on $ X $ and let
\begin{equation}
  \mathcal{I}_n(g) = \frac 1n \sum_{ j = 1 }^n g(\bx_j)
  .
\end{equation}
Then a simple calculation gives us
\begin{equation}
  \EE\big[ | I(g) - \mathcal{I}_n(g) |^2 \big] 
= \frac{\operatorname{Var}(g)}n
\qquad 
  \text{and}
\qquad
  \operatorname{Var}(g) 
  = \int_X | g(\bx) |^2 \, d\bx 
  - 
  \left[ \int_X g(\bx) \, d\bx \right]^2.
\label{MC-rate}
\end{equation}
The $O(1/\sqrt{n})$ rate is independent of $ d \in \N $.
To reduce the error to a tolerance of $ \eps > 0 $, 
the number of samples needed must be of order $ \eps^{ - 2 } $.
This is a situation with $\beta=2$ in \eqref{cost}.

The value of $\alpha$ is more subtle. We need to examine specific classes of examples in different
dimensions.  For example, one can ask about the value of $\alpha$ for the Ising or Heisenberg model in
statistical physics.  At the moment, results in this direction are still quite sparse.

Algorithms and results similar to those of the approximative computation of integrals 
in \eqref{eq:integral_intro}--\eqref{MC-rate} above 
have been developed in the case of linear PDEs of the Kolmogorov type but, 
for a very long time,  little progress was made in developing algorithms with quantitatively similar computational complexity 
for high dimensional nonlinear PDEs and this has impeded advances in several fields such as  
optimal control and quantum mechanics.
Things have changed dramatically 
in the last few years 
with the appearance of 
so-called {\it full history recursive 
multilevel Picard approximation} methods \cite{E2016multilevel,E2019multilevel,Hutzenthaleretal2018arXiv} 
and a host of machine learning-based algorithms for high dimensional PDEs 
beginning with the {\it Deep BSDE method} \cite{EHanJentzen2017,HanJentzenE2018}. 
Full history recursive multilevel Picard approximation methods 
(in the following we abbreviate 
\emph{full history recursive multilevel Picard} by MLP) 
are some recursive nonlinear variants of the
classical Monte-Carlo approximation methods 
and, in that sense, MLP approximation methods are 
nonlinear Monte-Carlo approximation methods. 
For every arbitrarily small $ \delta \in (0,\infty) $ it has been shown 
that MLP approximation methods achieve \eqref{cost} with 
$ \alpha = 1 $ and $ \beta = 2 +\delta $ for a wide class of
nonlinear (parabolic) equations (see Section~\ref{sec:MLP} below for details). 
Although a complete theory is still lacking, the Deep BSDE method has demonstrated very
robust performance in practice for a range of problems and has been extended in many different ways.
These developments will be reviewed below. 

Along with the work on developing algorithms, there has been some effort to develop an open-source
platform where codes, review papers, and other information can be shared.  The interested reader should consult:
\url{http://deeppde.org}.

Before launching into a more detailed discussion, 
it is useful to review briefly the two main ideas that 
we focus on in this paper:
machine learning approximation methods 
(see Section~\ref{sec:intro_deep_learning} below) 
and MLP approximation methods 
(see Section~\ref{sec:intro_MLP} below).

\subsection{A brief introduction of supervised learning}
\label{sec:intro_deep_learning}

The basic problem in supervised learning is as follows:
Given a natural number $ n \in \N $
and a sequence 
$ ( \bx_j, y_j ) = ( \bx_j , f^*( \bx_j ) ) $,
$ j \in \{ 1, 2, \dots, n \} $, 
of pairs of input-output data,  
we want to  recover the target function $ f^* $ as accurately as possible. 
We will assume that the input data $ \bx_j $, $ j \in \{ 1, 2, \dots, n \} $, 
is sampled from the probability distribution $ \mu $ on $ \R^d $.

{\bf Step 1. Choose a hypothesis space}. This is a set of trial functions $ \cH_m $
where $ m \in \N $ is a natural number that is strongly related to the dimensionality of $\cH_m$.
One might choose piecewise polynomials or wavelets.
In modern machine learning the most popular  choice is neural network functions.
Two-layer neural network functions (one input layer, one output layer that usually does not count, and one hidden layer) take the form 
\begin{equation}
f_m(\bx,  \theta) =
\frac 1 m \sum_{ j = 1 }^m a_j \sigma( \left< w_j , x \right> )
\end{equation}
where $ \sigma \colon \R \to \R $ is a fixed scalar nonlinear function 
and where $ \theta = ( a_j, w_j )_{ j \in \{ 1, 2, \dots, m \} } $ are the parameters to be optimized (or trained). 
A popular example for the nonlinear function $ \sigma \colon \R \to \R $ is the rectifier function (sometimes also referred to 
as ReLU (rectified linear unit) activation function in which case we have for all 
$ z \in \R $ that 
$ \sigma(z) = \max\{ z, 0 \} $. 
Roughly speaking, deep neural network (DNN) functions are obtained if one composes two-layer neural network functions several times.
One important class of DNN models are residual neural networks (ResNet). They closely resemble discretized 
ordinary differential equations and take the form 
\begin{align}
\label{ResNet}
z_{l+1} & = z_l + \frac 1 {LM} \sum_{j=1}^M a_{j,l}\sigma( \left< z_l, w_{j,l} \right> ), 
\qquad
z_0 = V{\bx}, 
\qquad
f_L(\bx, \theta) = \left< \alpha, z_L \right>
\end{align}
for $ l \in \{ 0, 1, \dots, L - 1 \} $ where $ L, M \in \N $. 
Here the parameters are 
$ 
  \theta 
  = ( \alpha, V, 
    ( a_{j,l} )_{ j, l } 
    , 
    ( w_{ j, l } )_{ j, l }
  ) 
$.
ResNets are the model of choice for truly deep neural network models.

{\bf Step 2. Choose a loss function}. The primary consideration for the
choice of the loss function is to fit the data.  Therefore one most obvious
choice is the $L^2$ loss:
\begin{equation}
\hat{\mathcal{R}}_n(f) = \frac1n \sum_{j=1}^n | f(\bx_j) - y_j |^2 = \frac 1n \sum_{j=1}^n | f(\bx_j) - f^*(\bx_j) |^2.
\label{empirical-risk}
\end{equation}
This is also called the ``empirical risk''.
Sometimes we also add regularization terms.

{\bf Step 3. Choose an optimization algorithm}.  
The most popular optimization algorithms in machine learning are different versions of
the gradient descent (GD) algorithm, or its stochastic analog, the
stochastic gradient descent (SGD) algorithm. Assume that the objective function we aim to minimize is of the form
\begin{equation}
F(\theta) = \E_{\xi \sim \nu}\big[ l(\theta, \xi) \big]
\label{expectation}
\end{equation}
($\nu$ could be an empirical distribution on a finite set). 
The simplest form of SGD iteration takes the form 
\begin{equation}
\theta_{k+1} = \theta_k - \eta \nabla l(\theta_k, \xi_k),
\end{equation}
for $ k \in \N_0 $ where $ \xi_k $, $ k \in \N_0 = \{ 0, 1, 2, \dots \} $, is a sequence of 
i.i.d.\ random variables sampled from the distribution $ \nu $ 
and $ \eta $ is the learning rate which might also change during the course of the iteration.
In contrast GD takes the form
\begin{equation}
\theta_{k+1} = \theta_k - \eta \nabla \E_{\xi \sim \nu}\big[ l(\theta_k, \xi) \big].
\end{equation}
Obviously this form of SGD can be adapted to loss functions of the form \eqref{empirical-risk} which can also be
regarded as an expectation.
The DNN-SGD paradigm is the heart of modern machine learning.

\paragraph{High dimensional stochastic control problems.}

One of the earliest applications of deep learning to problems in scientific computing is 
for the stochastic control problem \cite{HanE2016deepcontrol}.
This example was chosen because of its close resemblance to the DNN-SGD paradigm in deep learning.
From an abstract viewpoint, DNN can be viewed as a (discrete) dynamical system, of which  ResNet is a good example.
SGD is a natural consequence when applying GD to stochastic optimization problems, in which the objective function is an expectation.

Consider the stochastic dynamical system:
\begin{equation}
  s_{t +1} = s_t+b_t(s_t,a_t)+\xi_{t +1}.
  \label{scontrol}
  \end{equation}
  Here $s_t, a_t, \xi_t $ are respectively the  state, the control, and the noise at time $t$. 
  We assume that the objective function for the control problem is of the form:
  \begin{align}
  \min_{\{a_t\}_{t =0}^{T-1}}\mathbb{E}\Big[\sum_{t =0}^{T-1} c_t(s_t,a_t(s_t))+c_T(s_T)\mid s_0\Big],
  \label{sobjective}
\end{align}
where $T$ is the time horizon and  $a_t = a_t(s_t)$ are the feedback controls. 
One can see  the close analogy between the stochastic
control problem and the DNN-SGD paradigm:  \eqref{scontrol} plays the role of \eqref{ResNet} and \eqref{sobjective} is in the
form of a stochastic optimization problem. In this analogy, the role of the training data is played by the noise $\{\xi_t\}$.

\begin{figure}[!ht]
\centering
  \includegraphics[width=0.56\textwidth]{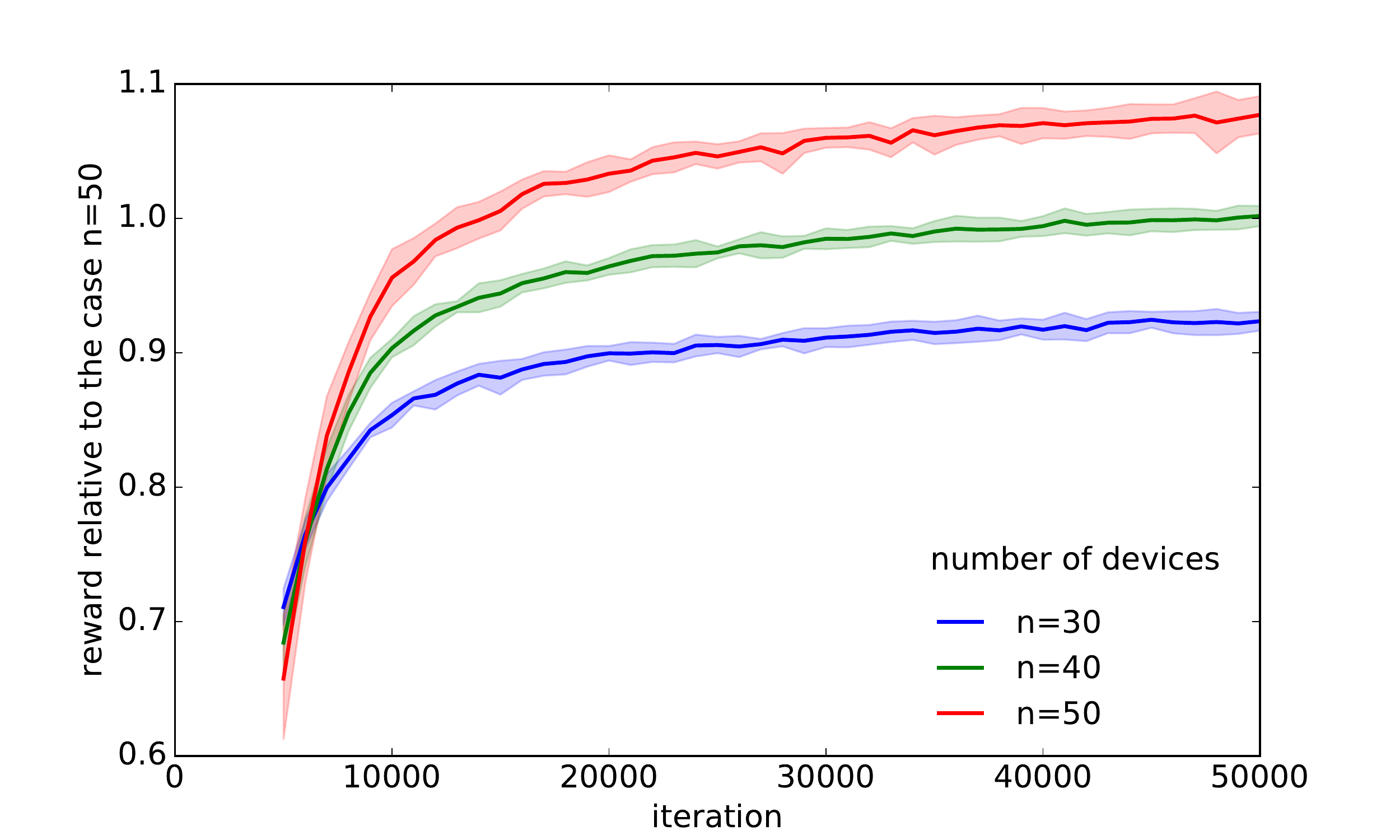}
  \captionof{figure}{Relative reward to the case of number of devices $n=50$ (with controls satisfying constraints strictly). The space of control function is $\mathbb{R}^{n+2}\rightarrow \mathbb{R}^{3n}$ for $ n \in \{ 30, 40, 50 \} $, 
  with multiple equality and inequality constrains. 
  The shaded area depicts the mean $\pm$ the standard deviation over five different random seeds. 
  Reprinted from \cite{HanE2016deepcontrol}.}
  \label{fig:energystorage}
\end{figure}

To develop an algorithm for this problem, one can  approximate the feedback control function $a_t$
by any machine learning model, in particular some neural network model:
  \begin{equation}
    a_t(s_t)\approx a_t(s_t|\theta_t),
  \end{equation}
  where $\theta_t$ is the parameter to be trained at time $t$.
  The loss function can be defined by
  \begin{equation}
 L(\{\theta_t\}) =  \mathbb{E}\!\left[\sum_{t =0}^{T-1}c_t(s_t,a_t(s_t|\theta_t))+c_T(s_T) \right],
  \end{equation}
  which can be optimized using SGD over different random samples of $ \xi_t $, $ t \in \{ 1, 2, \dots, T \} $.
An example of energy storage is shown in Figure \ref{fig:energystorage}.
It is an allocation
problem,  with the objective being optimizing total revenue from multiple storage devices and a renewable
wind energy source while satisfying stochastic demand.
More details of the problem can be found in \cite{HanE2016deepcontrol,jiang2015approximate}.

\subsection{A brief introduction to multilevel Picard approximation methods}
\label{sec:intro_MLP}

Despite the great performance of deep learning-based approximation
schemes in various numerical simulations, until today, 
there is no rigorous mathematical analysis in the scientific literature 
which proves (or disproves) 
the conjecture that there exists a deep learning-based approximation method 
which overcomes the curse of dimensionality 
in the numerical approximation of PDEs. 
However, for MLP approximation methods 
it has been proven in the scientific literature 
that such approximation methods do 
overcome the curse of dimensionality in the numerical approximation
of semilinear PDEs with general time horizons 
and, in particular, the convergence results for MLP approximation methods 
(see 
\cite{becker2020arxiv2005,E2016multilevel,Hutzenthaleretal2018arXiv,
hutzenthaler2019arxiv1903,
Becketal2019MLP_nonlip_arXiv,giles2019generalised,
beck2020arxiv2003,
hutzenthaler2019arxiv1912,hutzenthaler2020multilevel,
hutzenthaler2020lipschitz} and Section~\ref{sec:MLP} below for details) 
have revealed that semilinear PDEs with general time horizons 
can be approximated without the curse of dimensionality.

Let us briefly illustrate this in the case of semilinear heat PDEs 
with bounded initial value functions. 
Let $ T, c \in (0,\infty) $, 
let $ f \colon \R \to \R $ be Lipschitz continuous, 
for every $ d \in \N $ let 
$ u_d \in C^{ 1, 2 }( [0,T] \times \R^d, \R ) $, 
and assume for every 
$ d \in \N $, $ t \in [0,T] $,
$ x \in \R^d $
that 
$
  | u_d(t,x) | \leq c
$
and 
\begin{equation}
\label{eq:MLP_allen_cahn_pde_intro}
  (\tfrac{\partial}{\partial t}u_d)(t,x) 
  = 
  (\Delta_x u_d)(t,x) 
  + 
  f(u_d(t,x)) .
\end{equation} 
In the linear case where $ f \equiv 0 $ vanishes, 
it is known for a long time that 
classical Monte-Carlo approximation methods 
can approximate $ u_d( T, 0 ) \in \R $ 
with $ \alpha = 1 $ and $ \beta = 2 $ 
in \eqref{cost}. 
In the general nonlinear case, 
classical Monte-Carlo approximation methods 
are not applicable anymore 
but 
it has recently been shown in the scientific literature 
(see 
Hutzenthaler et al.~\cite{Hutzenthaleretal2018arXiv} 
and Theorem~\ref{thm:MLP_analysis} below) 
that 
for every arbitrarily small $ \delta \in (0,1] $ 
it holds that 
MLP approximation methods 
can approximate $ u_d( T, 0 ) \in \R $ 
in the general nonlinear case 
with $ \alpha = 1 $ and $ \beta = 2 + \delta $ 
in \eqref{cost}. 
The convergence results for MLP approximation methods 
in the scientific literature have thus revealed that 
semilinear heat PDEs can, up to an arbitrarily small 
polynomial order of convergence, been solved 
with the same computational complexity 
as linear heat PDEs.

In the following we briefly sketch some of the main ideas 
in the derivation of MLP approximation schemes. 
One of the key ideas in the derivation and the convergence analysis 
of the MLP approximation scheme is to rewrite 
the PDE in \eqref{eq:MLP_allen_cahn_pde_intro} as a 
stochastic fixed point equation. More formally, we note that 
\eqref{eq:MLP_allen_cahn_pde_intro} ensures that 
for all $ t \in [0,T] $, $ x \in \R^d $ it holds that
\begin{equation}
\label{eq:fixed_point_intro_1}
  u(t,x) 
  =
  \E\!\left[ 
    u\big( 
      0, x + \sqrt{2} W_t 
    \big) 
    +
    \int_0^t 
    f\big( 
      u(s, x + \sqrt{2} W_{ t - s } ) 
    \big)
    \, ds
  \right]
  .
\end{equation}
where $ W \colon [0,T] \times \Omega \to \R^d $ 
is a standard Brownian motion 
on a probability space $ ( \Omega, \mathcal{F}, \P ) $ 
(cf., e.g., Beck et al.~\cite[Theorem~1.1]{beck2020nonlinear}). 
Now we can also write \eqref{eq:fixed_point_intro_1} 
as the fixed point equation 
$ u = \Phi( u ) $
where $ \Phi $ is the self-mapping on the set of all bounded functions 
in $ C( [0,T] \times \R^d, \R ) $ which is described through 
the right hand side of \eqref{eq:fixed_point_intro_1}. 
Using $ \Phi $ one can define Picard iterates $ \mathfrak{u}_n $, $ n \in \N_0 $, through the recursion 
that for all $ n \in \N $ it holds that
$ \mathfrak{u}_0 = 0 $ and $ \mathfrak{u}_n = \Phi( \mathfrak{u}_{ n - 1 } ) $. 
In the next step we note that a telescoping sum argument shows that for all $ n \in \N $ it holds that
\begin{equation}
\label{eq;intro_telescopic}
  \mathfrak{u}_n 
  = \mathfrak{u}_1 + \sum_{ k = 1 }^{ n - 1 } \left[ \mathfrak{u}_{ k + 1 } - \mathfrak{u}_k \right]
  = \Phi( 0 ) + \sum_{ l = 1 }^{ n - 1 } \left[ \Phi( \mathfrak{u}_k ) - \Phi( \mathfrak{u}_{ k - 1 } ) \right] 
  .
\end{equation}
MLP approximations are then derived by approximating 
the expectations in \eqref{eq;intro_telescopic} 
within the fixed point mapping $ \Phi $ by means of Monte-Carlo approximations 
with different levels of accuracy.

The procedure in \eqref{eq;intro_telescopic} is inspired by multilevel Monte Carlo (MLMC) approximations 
in the scientific literature (see Heinrich~\cite{heinrich1998monte}, Heinrich \& Sindambiwe~\cite{heinrich1999monte}, Giles~\cite{giles2008multilevel} and, e.g.,~\cite{heinrich2001multilevel,giles2015multilevel} and the references mentioned therein). 
There are, however, also several differences when comparing MLP approximations 
to ``classical'' MLMC approximations. In particular, we note that MLP approximations are 
full history recursive in the sense that for all $ n \in \N $ we have that computations of realizations of MLP approximations in the $ n $-th iterate require realizations for MLP approximations 
in the 1st, 2nd, $ \dots $, $ (n - 1) $-th iterate 
the analysis of MLP approximations (see \eqref{eq:mlp} 
in Theorem~\ref{thm:MLP_analysis} below for details). 
Taking this into account, the convergence analysis for MLP approximations in the scientific literature 
turns out to be more much subtle, and we refer to Section~\ref{sec:MLP} below 
for a sketch of the some of the main ideas for the complexity analysis of MLP approximations
and also for references to research articles on MLP approximations in the scientific literature.

In the comparison between 
classical linear Monte-Carlo methods 
and MLP approximation methods 
in \eqref{eq:MLP_allen_cahn_pde_intro} above 
we have restricted ourselves just for simplicity 
to the problem of approximating 
semilinear \emph{heat} PDEs with \emph{bounded} solutions 
at the space-time point $ t = T $, $ x = 0 $ 
and similar results hold in the case of much more general classes 
of PDEs and more general approximation problems. 
Until today, MLP approximation schemes are the only approximation schemes 
in the scientific literature for which it has been proved
that they do overcome the curse of dimensionality 
in the numerical approximation 
of semilinear PDEs with general time horizons. 
We refer to Section~\ref{sec:MLP} for further details and, in particular, 
for a comprehensive literature overview on research articles 
on MLP approximation methods.

\section{General remarks about algorithms for solving PDEs in high dimensions}
\label{sec:general_remarks}

We begin with a brief overview of the algorithms that have been developed for high dimensional PDEs.

\vspace{.1in}
\noindent
{\bf Special classes of PDEs:} 
The representative special classes of PDEs 
which we review within this subsection  
are second-order linear parabolic PDEs of the Kolmogorov type 
and first-order Hamilton-Jacobi equations. 

Consider the linear parabolic PDE with a terminal condition specified at $ t = T $ described by
  \begin{align}
  \label{linear-PDE}
    \frac{ \partial u}{ \partial t } 
    + \frac{1}{2} \! \Tr\!\big( \sigma\sigma^{\operatorname{T}}( \operatorname{Hess}_x u) \big)
    + \left< \nabla_x u, \mu \right> + f = 0, \quad u(T, \cdot) = g(\cdot)
  \end{align}
Our objective is to compute $ u(0, \cdot)$. 
For this purpose, we consider the diffusion process described by 
the stochastic differential equation (SDE)
\begin{equation}
  d\bX_t = \mu( t, \bX_t )\, dt + \sigma( t, \bX_t ) \, W_t .
\end{equation}
The solution to the PDE in \eqref{linear-PDE} can be expressed 
as an expectation in the sense that
\begin{equation}
  u(0,\bx) = \E\!\left[ g(\bX_T)+\int_0^T f(s, \bX_s) ds \, \bigg| \bX_0=\bx \right].
\end{equation}
This is the classical  Feynman-Kac formula~\cite{Karatzas1998,oksendal2013stochastic}.
Using this formula, one can readily evaluate $u(0, \cdot)$ using Monte Carlo without suffering from CoD.

In a similar spirit, solutions of the Hamilton-Jacobi equations 
can also be expressed using the Hopf formula.  
Consider the PDE
\begin{equation}
  \frac{ \partial u }{ \partial t } + H( \nabla u ) = 0, \qquad u(\bx, 0) = g(\bx).
\end{equation}
Assume that $H$ is convex.  Then we have the Hopf formula:
\begin{equation}
u(\bx, t) = \inf_{y} \left( g(y) + t H^*\Big(\frac{\bx - y}{t}\Big) \right).
\end{equation}
where $H^*$ is the Legendre transform of $ H $ (see Evans~\cite{Evans2010partial}). 
The right hand side of the above equation can be computed using
some optimization algorithms. 
A particularly attractive algorithm along these lines 
was developed in Darbon \& Osher~\cite{DarbonOsher2016}.

\vspace{.1in}
\noindent
{\bf Control problems:}
      The first application of deep learning to solving scientific computing problems in high dimensions was in the area of
      stochastic control.  In 2016, Han and E \cite{HanE2016deepcontrol} developed  a deep neural network-based algorithm for stochastic control problems.
      The reason for choosing stochastic control  was its very close analogy with the setup of deep learning, 
      as we will see later (see Section~\ref{sec:control} below).
       Deep learning-based algorithms for deterministic control problems were first developed in
       \cite{nakamura2019adaptive}.

\vspace{.1in}
\noindent
{\bf Schr\"odinger equation for spins and electrons:}
  In an influential paper, Carleo and Troyer  developed an algorithm for solving the Schr\"odinger
  equation for spins using the restricted Boltzmann machine (RBM) as the trial function.
  The variational Monte Carlo (VMC) approach was used for ground-state calculations.  To solve the dynamic equation,
  the  least square approach was used, i.e., the total integral of the square of the residual was used as the loss function \cite{carleo2017solving}.
  
  For many-electron Schr\"odinger equation, the story is quite different. The configuration space is now continuous, instead of being discrete.
  In addition, the wave function should satisfy the anti-symmetry constraint. This has proven to be a difficult issue in solving the 
  Schr\"odinger equation.
  In \cite{han2019solving}, Han, Zhang, and E  developed a deep neural network-based methodology for computing the ground states.
  Compared with traditional VMC, a permutation-symmetric neural network-based ansatz is used for the Jastrow factor. The
  anti-symmetric part was treated in a rather simplified fashion. This has been improved in the work~\cite{luo2019backflow, pfau2019ab, Hermann2019deep}. 
   
  Despite these progresses, it is fair to say that we are still at an early
  stage for developing  neural network-based algorithms for the many-body Schr\"odinger equation.
  Since the issues for solving the Schr\"odinger equation are quite specialized,  we will not discuss them further in this review.
    
\vspace{.1in}
\noindent
{\bf  Nonlinear parabolic PDEs: }
The first class of algorithms developed for general nonlinear parabolic PDEs with general time horizons 
in really high dimensions ($d\geq 40$) is the 
{\it multilevel Picard method} \cite{E2016multilevel,E2019multilevel,Hutzenthaleretal2018arXiv}. 
At this moment, this is also the only algorithm for which a relatively complete theory has been established 
(see Section~\ref{sec:MLP} below).
Among other things, this theory offers a proof that the MLP method overcomes CoD.
Shortly after, E, Han, and Jentzen developed 
the deep neural network-based {\it Deep BSDE method},
by making use of the connection between nonlinear parabolic equations and backward stochastic differential
equations (BSDE) \cite{EHanJentzen2017, HanJentzenE2018}. 
This was the first systematic application of deep learning to solving general high dimensional PDEs. 
Later, Sirignano and Spiliopoulos developed an alternative deep learning-based algorithm using the least squares approach 
\cite{Sirignano2018dgm}, extending the work of Carleo and Troyer  \cite{carleo2017solving} to general PDEs. 
Such deep learning-based approximation methods for PDEs have also been 
extended in different ways and to other parabolic and even elliptic problems; 
see, e.g., \cite{beck2019machine,BeckerCheridito2019,raissi2018forward,BeckBeckerCheridito2019,BeckerCheriditoJentzen2019,HurePhamWarin2019,ChanMikaelWarin2019,Henry2017deep,han2019solving,ji2020three,han2020solving}.

Some special semilinear parabolic PDEs can be formulated in terms of branching processes. 
One such example is the Fisher-KPP (Kolmogorov-Petrovski-Piscounov) equation \cite{Fisher1937wave,KPP1937,Mckean1975application}. 
For such PDEs, Monte Carlo methods can be developed, and such Monte Carlo approximation algorithms overcome the CoD 
(see \cite{skorokhod1964branching,watanabe1965branching,Henry-Labordere2012,
Henry-Labordere2014,chang2016branching,warin2017variations,henry2019branching}) 
in the specific situation where the time horizon and/or the PDE nonlinearity is sufficiently small.

\vspace{.1in}   
\noindent
{\bf Variational problems:}
It is fairly straightforward to construct neural network-based algorithms for solving variational problems.  One way of doing this is
simply to use the Ritz formulation. The {``Deep Ritz} method'', to be discussed below, is such an example; 
see E \& Yu~\cite{E2018deep} and also Khoo et al.\ \cite{Khoo2019solving}.
It is natural to ask whether one can develop a similar Galerkin formulation,  i.e., 
using a weak form of the PDE. In fact, \cite{E2018deep} was written
as a preparation for developing what would be a ``Deep Galerkin method''. 
However, formulating robust neural network algorithms using
a weak form has proven to be quite problematic.  The difficulty seems to be analogous to the ones encountered in 
generative adversarial networks (GAN); cf.\ Goodfellow et al.\ \cite{Goodfellow2014generative}. 
For some progress in this direction we refer to the article Zhang et al.\ \cite{Zang2020weak}.
It should also be noted that even though the deep learning-based methodology proposed in 
Sirignano \& Spiliopoulos~\cite{Sirignano2018dgm} was named as a 
``Deep Galerkin method'', the methodology in \cite{Sirignano2018dgm} 
is somehow based on a least square formulation rather than a Galerkin formulation. 

\vspace{.1in}
\noindent
{\bf Parametric PDEs:}
One of the earliest applications of deep learning to PDEs is in the study of parametric PDEs.
In \cite{Khoo2020solving}, Khoo, Lu, and Ying developed a methodology for solving low dimensional
PDEs with random coefficients in which the neural network models are used to parametrize the random coefficients. Recently the neural networks are also applied to solve low-dimensional stochastic PDEs~\cite{zhang2020learning}.
This is a promising direction thought it will not be covered in this review.
Another closely related area is solving inverse problems governed by PDEs, which is intrinsically high dimensional as well. Recent works~\cite{raissi2019physics,fan2019solving,fan2020solving,khoo2019switchnet,chen2020physics} have demonstrated the advantages of approximating the forward and inverse maps with carefully designed neural networks.

\vspace{.1in}
\noindent
{\bf Game theory}
A stochastic game describes the behavior of a population of interactive agents among which everyone makes his/her optimal decision in a common environment. Many scenarios in finance, economics, management science, and engineering can be formulated as stochastic games. 
With a finite number of agents, the Markovian Nash equilibrium of a game is determined by a coupled system of parabolic PDEs. 
To solve these problems, Han et al. extend the Deep BSDE method in \cite{Han2019deep,han2020convergence} with the idea of fictitious play. 
In a different direction, with an infinite number of agents and no common noise, one can use the mean-field game theory developed in \cite{LaLi1:2006,LaLi2:2006,LaLi:2007,HuMaCa:06,HuCaMa:07} to reduce the characterization of the Nash equilibrium to two coupled equations: a 
Hamilton-Jacobi-Bellman equation and a Fokker-Planck equation. Neural network-based algorithms have been developed in \cite{Carmona2019convergence1,Carmona2019convergence2,Ruthotto2020machine,Lin2020apac} to solve these equations.

\vspace{.1in}
Besides the literature mentioned above, certain deep learning-based approximation methods for PDEs have been proposed (see, e.g.,~\cite{berg2018unified,dockhorn2019discussion,farahmand2017deep,fujii2019asymptotic,goudenege2019variance,jacquier2019deep,lye2020deep,magill2018neural,nusken2020solving,pham2019neural,raissi2018deep,cai2018approximating,zhang2019quantifying}) and various numerical simulations for such methods suggest that deep neural network approximations might have the capacity to indeed solve high dimensional problems efficiently. Actually, the attempts of applying neural networks for solving PDEs can be dated back to the 90s (cf., e.g.,~\cite{lee1990neural,uchiyama1993solving,lagaris1998artificial,jianyu2003numerical}), nevertheless, with a focus on low-dimensional PDEs.
Apart from neural networks, there are also other attempts in literature in solving high dimensional PDEs with limited success (see, e.g., \cite{bally2003quantization,von2004numerical,zhang2004numerical,bouchard2004discrete,gobet2005regression,delarue2006forward,bender2007forward,gobet2008numerical,briand2014simulation,geiss2016simulation,gobet2010solving,fahim2011probabilistic,crisan2010probabilistic,crisan2010monte,crisan2012solving,crisan2014second,labart2013parallel,pham2015feynman,guo2015monotone,gobet2016stratified,gobet2016approximation,gobet2016linear,warin2018nesting,warin2018monte,ruszczynski2017dual,billaud2018stochastic}).

This review will be focused on nonlinear parabolic PDEs and related control problems.
There are two main reasons for choosing these topics. The first is that these classes of problems are fairly general
and have general interest.  The second is that the study of high dimensional problems  is in better shape for these classes of problems, compared to others (e.g., the Schr\"odinger equation discussed above).

We should also mention that the heart of reinforcement learning is solving approximately the Bellman equation, even though reinforcement learning
algorithms are not always formulated that way. The dimensionality in these problems is often very high. This is another topic that will not be
covered in this review.

\section{The Deep BSDE method}
The Deep BSDE method was the first deep learning-based  numerical algorithm  for solving general nonlinear parabolic PDEs in high dimensions~\cite{EHanJentzen2017, HanJentzenE2018}.
It begins by reformulating the PDE as a stochastic optimization problem.  This is done with the help of BSDEs, hence the name ``Deep BSDE method''.   As a by-product, the Deep BSDE method is also an efficient algorithm
for solving high dimensional BSDEs.

\label{sec:nonlinear_PDEs}
\subsection{PDEs and BSDEs}
Consider the semilinear parabolic PDE
\begin{equation}
    \label{eq:PDE}
    \frac{ \partial u}{ \partial t } 
    + \frac{1}{2} \! \Tr\!\big( \sigma \sigma^{\operatorname{T}} (\operatorname{Hess}_{\bx} u) \big)
    +\left< \nabla u, \mu \right>
    +f\big( t, \bx, u, \sigma^{\operatorname{T}} (\nabla_x u) \big) = 0, \quad u(T, \bx) = g(\bx).
\end{equation}
In the same way as in Section~\ref{sec:general_remarks} above, we consider the diffusion process
\begin{equation}
\bX_t = \xi + \int_0^t\mu( s, \bX_s )\, ds +\int_0^t \sigma( s, \bX_s ) \, d W_s.
\end{equation}
Using It\^{o}'s lemma, we obtain that
\begin{equation}
\begin{split}
&
  u(t, \bX_t) - u(0, \bX_0)
\\
= &~-\int_0^t f\big( 
  s, \bX_s, u(s,\bX_s), [ \sigma( s, \bX_s ) ]^{ \operatorname{T} } ( \nabla_x u )( s, \bX_s )
  \big) \, ds \\
  &~+ \int_0^t [ \nabla u( s, \bX_s ) ]^{ \operatorname{T} } \,\sigma( s, \bX_s )\, d W_s.
  \end{split}
\end{equation}
To proceed further, we recall the notion of 
backward stochastic differential equations (BSDEs) 
\begin{empheq}[left=\empheqlbrace]{align}
\label{eq:BSDE_1}
    &\bX_t = \xi + \int_0^t\mu( s, \bX_s )\, ds +\int_0^t \sigma( s, X_s ) \, dW_s, \\
    &Y_t = g( X_T ) + \int_t^T f( s, X_s, Y_s, Z_s ) \, ds 
    - \int_t^T ( Z_s )^{ \operatorname{T} } \, dW_s
\label{eq:BSDE_2}
\end{empheq}
introduced by Pardoux and Peng \cite{Pardoux1992}.
It was shown in \cite{Pardoux1992,Pardoux1999} that there is 
an up-to-equivalence unique 
adapted stochastic process $ ( X_t, Y_t, Z_t ) $, $ t \in [0,T] $, 
with values in
$ \R^d \times \R \times \R^d $ that satisfies the pair of stochastic equations 
in \eqref{eq:BSDE_1}--\eqref{eq:BSDE_2} above. 

The connection between the BSDE in \eqref{eq:BSDE_1}--\eqref{eq:BSDE_2} 
and the PDE in \eqref{eq:PDE} is as follows \cite{Pardoux1992,Pardoux1999}.
Let $ u \colon [0,T] \times \R^d \to \R $ be a solution of the PDE in \eqref{eq:PDE}. 
If we define
\begin{equation}
\label{eq:nonlinear_Feynman_Kac}
  Y_t = u( t, X_t )
\qquad  
  \text{and}
\qquad 
  Z_t = [ \sigma( t, X_t ) ]^{ \operatorname{T} }  
  ( \nabla_x u )( t, X_t ).
\end{equation}
Then $ ( Y_t, Z_t ) $, $ t \in [0,T] $, is a solution for the BSDE in \eqref{eq:BSDE_1}--\eqref{eq:BSDE_2}.
With this connection in mind, one can formulate the PDE problem as 
the following variational problem:
\begin{align}
&\inf_{Y_0,\{Z_t\}_{0\le t \le T}} \E\big[ |g(X_T) - Y_T|^2 \big], \\
&s.t.\quad X_t = \xi + \int_{0}^{t}\mu(s,X_s)\, \,ds + \int_{0}^{t}\Sigma(s,X_s)\, dW_s, \\
&\hphantom{s.t.}\quad Y_t = Y_0 - \int_{0}^{t}h(s,X_s,Y_s,Z_s)\,  ds + \int_{0}^{t}(Z_s)\transpose\, dW_s.
\end{align}
The  minimizer of this variational problem is the solution to the PDE and vice versa.

\subsection{The Deep BSDE Method}

A key idea of the Deep BSDE method 
is to approximate the unknown functions
    $X_0 \mapsto u(0, X_0)$ {and}
    $X_t \mapsto [ \sigma(t,X_t) ]^{ \operatorname{T} } ( \nabla_x u )(t,X_t)$
    by feedforward neural networks $\psi$ and $\phi$.  To that purpose, we 
    work with the variational formulation described above and  discretize time, say using the Euler 
scheme on a grid $0 = t_0 < t_1 < \ldots < t_N = T$:
    \begin{align}
    &\inf_{\psi_0, \{\phi_n\}_{n=0}^{N-1}} \E |g(X_T) - Y_T|^2, \\
    &s.t.\quad X_0 = \xi, \quad Y_0 = \psi_0(\xi), \\
    &\hphantom{s.t.}\quad X_{t_{n+1}} = X_{t_i} + \mu(t_n,X_{t_n})\Delta t + \sigma(t_n,X_{t_n})\Delta W_n, \\
    &\hphantom{s.t.}\quad Z_{t_n} = \phi_n(X_{t_n}), \\
    &\hphantom{s.t.}\quad Y_{t_{n+1}} = Y_{t_n} - f(t_n,X_{t_n},Y_{t_n},Z_{t_n})\Delta t + (Z_{t_n})\transpose\Delta W_n.
    \end{align}
 At each time slide $t_n$, we associate a subnetwork.
 We can stack all these subnetworks together to form a deep composite neural 
    network. This network takes the paths 
    $\{ X_{ t_n } \}_{ 0 \leq n \leq N }$ and 
    $\{ W_{ t_n } \}_{ 0 \leq n \leq N }$ 
    as the input data and gives the final output, denoted by
        $\hat{u}( 
        \{ { X_{ t_n } } \}_{ 0 \leq n \leq N } , 
        \{ W_{ t_n } \}_{ 0 \leq n \leq N } )$, 
    as an approximation to 
    $u( t_N, X_{ t_N } )$. 
    
The error in the {matching of the given terminal condition} defines the  loss function 
      \begin{equation}
        l(\theta) = 
        \E\Big[
          \big|g( X_{ t_N } ) - \hat{u}\big(\{ X_{ t_n } \}_{ 0 \leq n \leq N } , \{ W_{ t_n } \}_{ 0 \leq n \leq N }\big)\big|^2
        \Big].
      \end{equation}

\begin{figure}[H]
\centering
\includegraphics[width=0.9\textwidth]{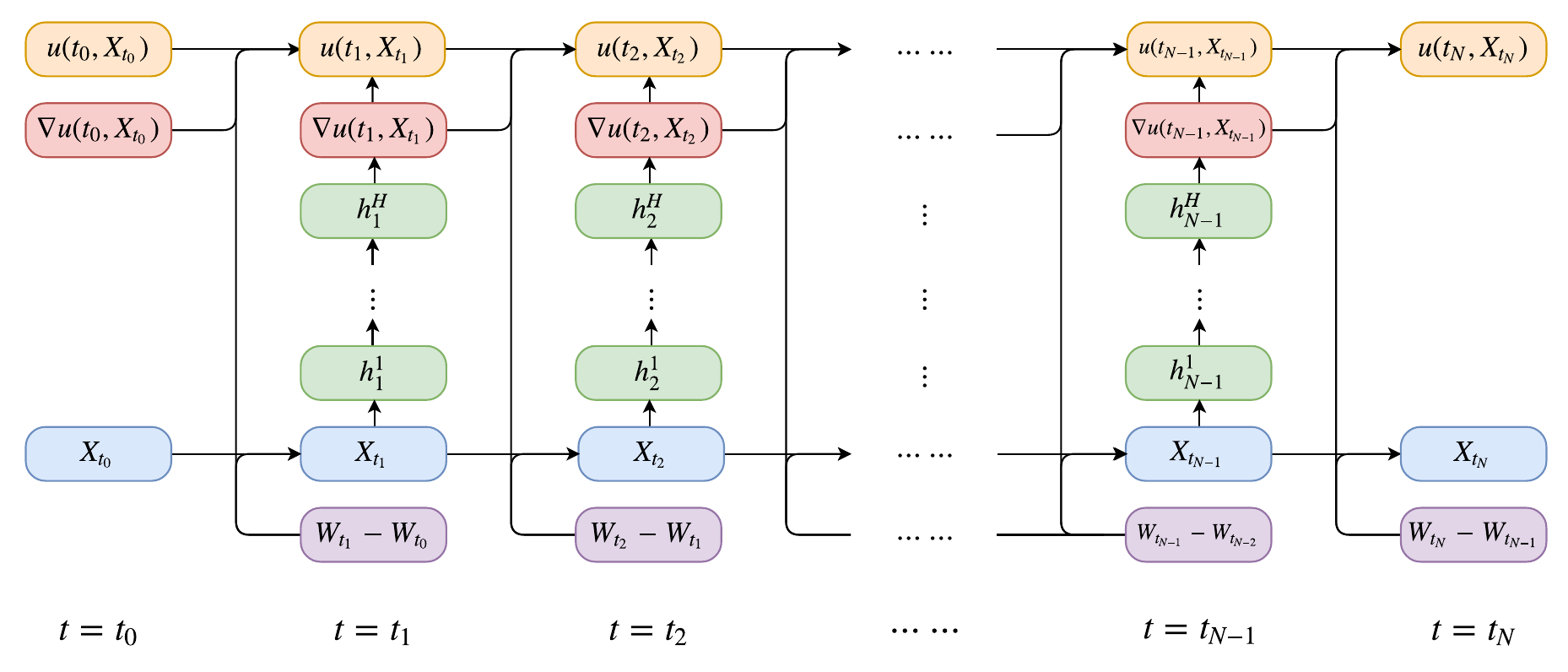}
\caption{Network architecture for solving parabolic PDEs. Each column corresponds to a subnetwork at time $t=t_n$. The whole network has $(H+1)(N-1)$ layers in total that involve free parameters to be optimized simultaneously. 
Reprinted from \cite{HanJentzenE2018}.}
\end{figure}

From the viewpoint of machine learning, this  neural network model has several interesting features.
\begin{enumerate}
\item It does not require us to generate training data beforehand.  The paths  $ \{ W_{ t_n } \}_{ 0 \leq n \leq N } )$ play
the role of the data and they are generated on the fly.  For this reason, one can think of this model as a model with an infinite amount of data.
\item For the same reason,  it is very natural to use stochastic gradient descent (SGD) to train the network.
\item The network has a very natural ``residual neural network'' structure embedded in the stochastic difference equations.
For example:
\begin{equation}
\begin{split}
  & u( t_{ n + 1 },X_{t_{n+1}}) - u( t_{ n },X_{t_{n}}) \\
\approx  & - f\big(t_n, X_{ t_n }, u( t_{ n },X_{t_{n}}),
    \phi_n(X_{t_n})
  \big)\,\Delta t_n+ (\phi_n(X_{t_n}))\transpose
    \Delta W_n.
\end{split}
\end{equation}
\end{enumerate}

\subsection{Some numerical examples}

Next we examine the effectiveness of the algorithms described above.  We will discuss two examples:
The first is a canonical benchmark problem, the linear-quadratic control problem (LQG).  The second is a nonlinear Black-Scholes model.
We use the simplest implementation of Deep BSDE:  Each subnetwork has $ 3 $ layers, with $ 1 $ input layer ($ d $-dimensional), $ 2 $ hidden layers (both $ d+10 $-dimensional), and $d$-dimensional output.
We choose the rectifier function (ReLU) as the activation function and optimize with the Adam method \cite{Kingma2015Adam}.
We will report the mean and the standard deviation of the relative error from  5 independent runs  with different random seeds.

\vspace{.1in}
\noindent
{\bf LQG (linear quadratic Gaussian)}

\begin{figure}[H]
    \centering
    \includegraphics[width=0.5\textwidth]{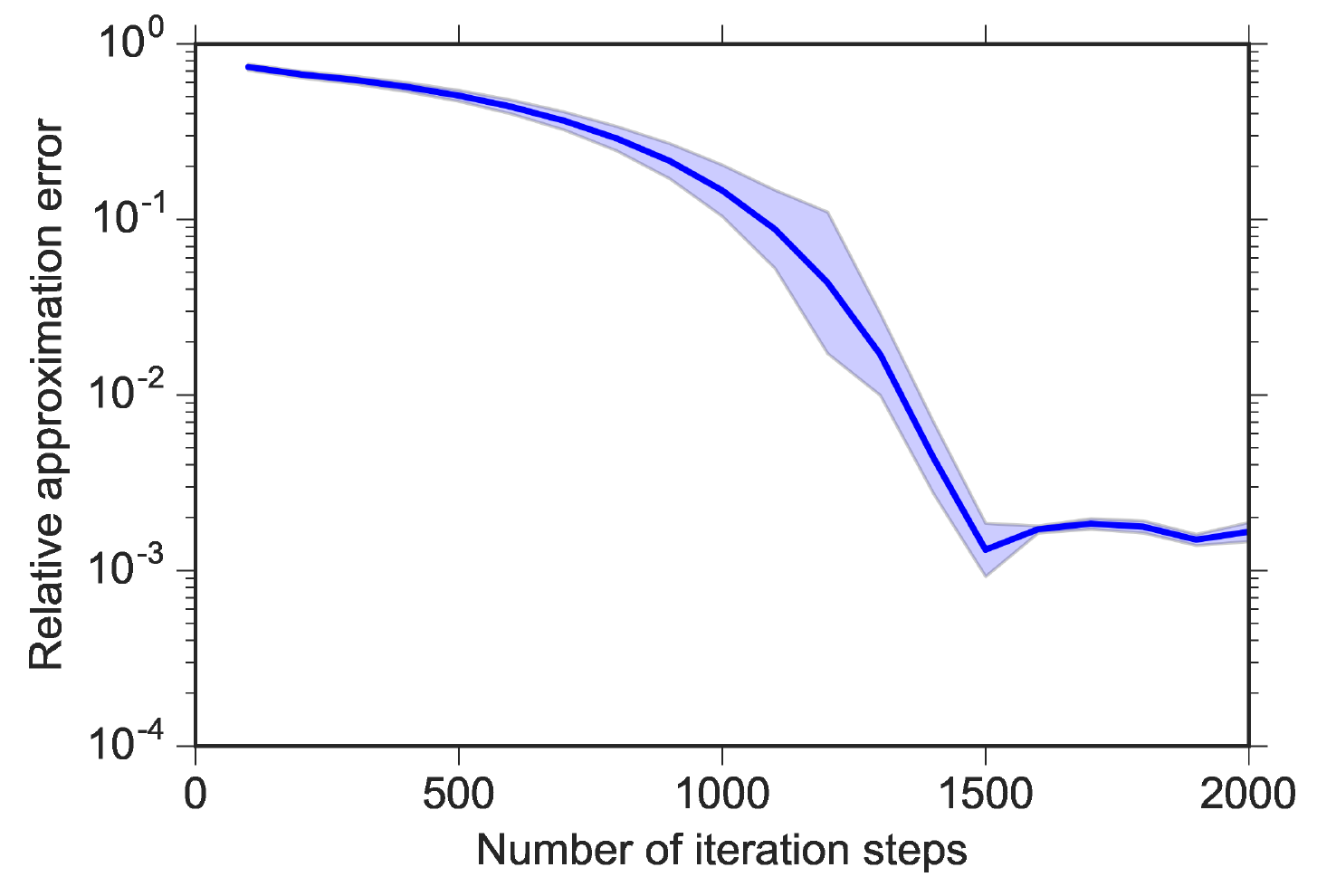}
    \includegraphics[width=0.47\textwidth]{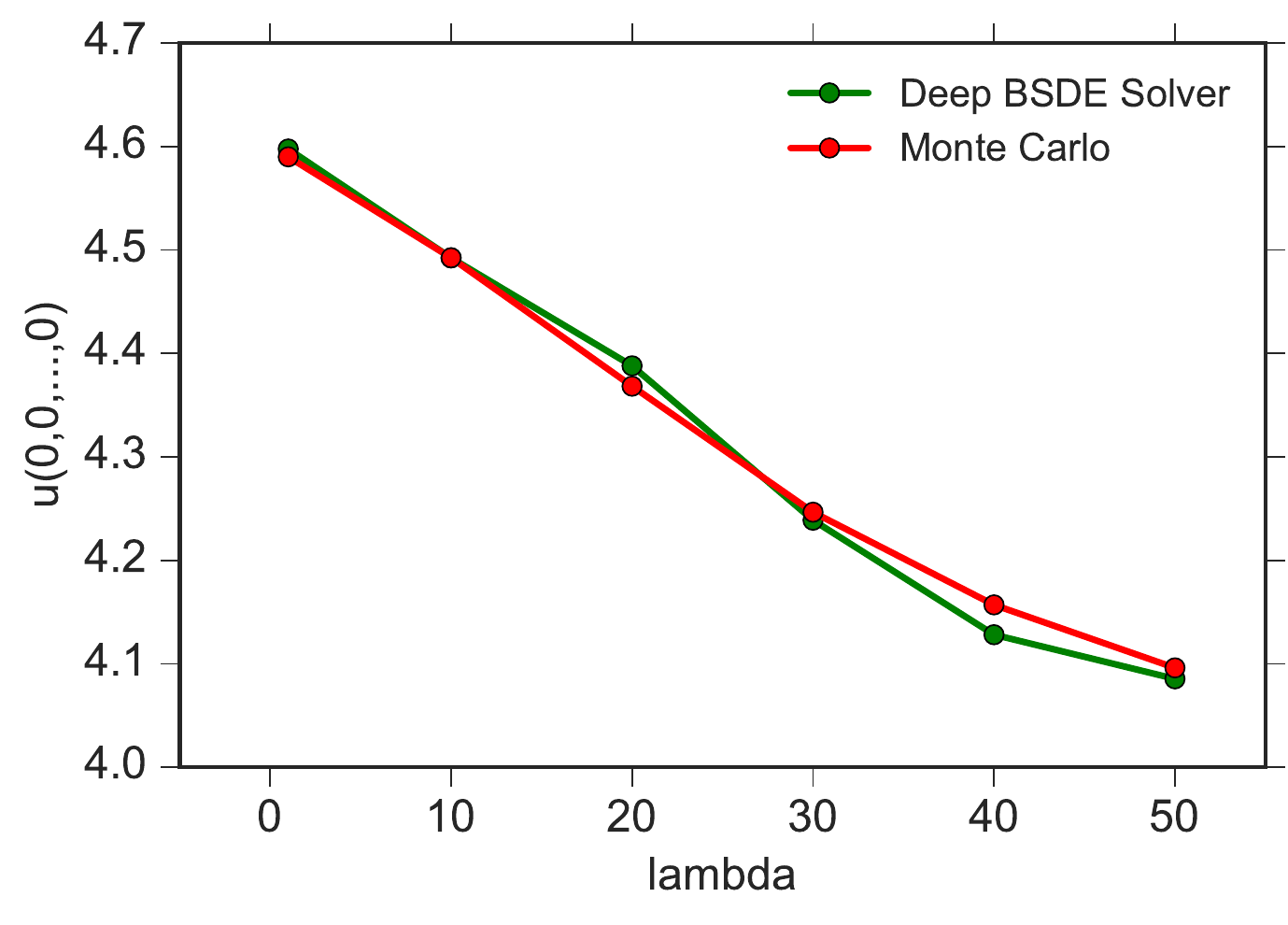}
    \caption{Left: Relative error of the Deep BSDE method for $ u( t{=}0, x{=}(0,\dots,0) )$ when $ \lambda = 1 $,
    which achieves $ 0.17\% $ relative error in a run time of 330 seconds. 
    Right: Optimal cost $u(t{=}0,x{=}(0,\dots,0))$ for different values of $\lambda$.
    The shaded area depicts the mean $\pm$ the standard deviation over five different random seeds. 
    Reprinted from \cite{HanJentzenE2018}.
    }
    \label{fig:LQR}
\end{figure}

Consider the stochastic dynamic model in 100 dimension:
\begin{equation}
            dX_t = 2\sqrt{\lambda}\,m_t\,dt+\sqrt{2}\,dW_t,
 \end{equation}
 with cost functional:
\begin{equation}
J( \{ m_t \}_{ 0 \leq t \leq T } ) =
          \E\big[
            \int_0^T \|m_t\|_2^2 \, dt + g(X_T)
          \big].
\end{equation}
        The associated HJB equation is given by
        \begin{equation}
     \frac{ \partial u}{ \partial t } + \Delta u  - \lambda \|\nabla u \|_2^2 = 0
        \end{equation}
        The solution to this HJB equation can be expressed as
\begin{equation}
                  u(t,x) = - \frac{ 1}{ \lambda }
                  \ln\!\bigg(
                    \E\Big[
                      \exp\!\Big(
                         - \lambda g( x + \sqrt{ 2 }W_{ T - t }  )
                      \Big)
                    \Big]
                  \bigg).
                \end{equation}
                This formula can be evaluated directly using Monte Carlo. Therefore this problem serves as a good model
                for validating algorithms. The results from the Deep BSDE method is shown in Figure \ref{fig:LQR}.
We see that the accuracy of the trained solution improves along the training curve before it saturates.

\vspace{.1in}
\noindent
{\bf Black-Scholes equation with default risk}

 The pricing model for financial derivatives  should take into account the whole basket of the underlies, which results in high dimensional PDEs.  In addition, 
the classical Black-Scholes model can and should be augmented by some important factors in real markets, including the 
 effect of default, transactions costs, uncertainties in the model parameters, etc.
 Taking into account these effects leads to nonlinear Black-Scholes type of models.

\begin{figure}[H]
  \centering
  \includegraphics[width=0.5\textwidth]{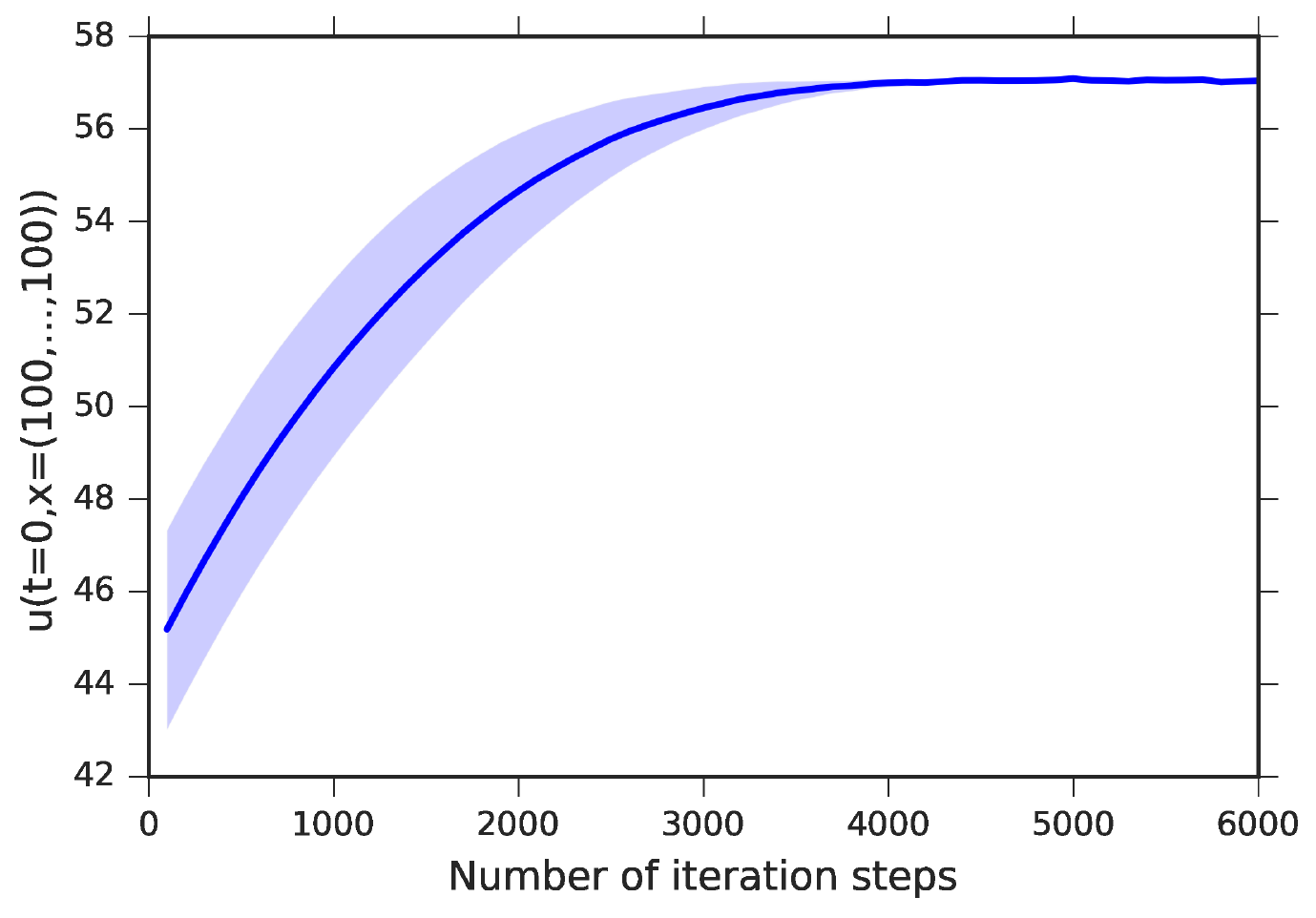}
  \caption{Approximation of $u(t{=}0,x{=}(100,\dots,100))$ as a function of the  number of iteration steps.
  The Deep BSDE method achieves a relative error of size $ 0.46\% $ in a runtime of $ 617 $ seconds.
  The shaded area depicts the mean $\pm$ the standard deviation over five different random seeds. 
  Reprinted from \cite{HanJentzenE2018}.
  }
  \label{fig:BlackScholes}
  \end{figure}

  We study a particular case of the recursive valuation model with default risk \cite{Duffie1996, Bender2017}. 
 The underlying asset price moves as a geometric Brownian motion, and the possible default is modeled by the first jump time of a Poisson process.
 The claim value is modeled by the nonlinear Black-Scholes model with
    \begin{equation}
      f\big( t, x, u(t,x), \sigma^{\operatorname{T}}( t, x ) \nabla u( t, x ) \big)
      = - \left( 1 - \delta \right) Q( u(t,x) ) \, u(t,x) - R \, u(t,x).
    \end{equation}
where $Q$ is some nonlinear function. %
We will consider the fair price of an European claim based on 100 underlying assets
    conditional on no default having occurred yet.  This leads to a problem with $d=100$.
Figure \ref{fig:BlackScholes} presents the results of Deep BSDE and multilevel Picard for this nonlinear Black-Scholes equation
for $d=100$.
Reported in the figure is the approximate solution at $t=0, $ $x=(100,\dots, 100)$.
For this problem we cannot find the ``exact solution''. Therefore we use the results of the two different methods to
calibrate each other.
  
\subsection{Analysis of the Deep BSDE method}
There is not yet a complete theory for the analysis of the Deep BSDE method.  We will review the existing results
that have been obtained so far.
Here instead of bounding the cost required  for reducing the error to certain tolerance $\eps$, we bound the error 
associated with certain hyper-parameters, such as the time step size $\Delta t$ and the size of the neural network models.
The basic strategy is to reduce the problem to bounding the generalization error for supervised learning
\cite{EMaWu2019machine}. In order to do that, we need to do the following: (1) Estimate the error in the time discretization. 
(2) Prove that the functions that need to be approximated using neural networks belong to the right function class
and bound their norms in that function class.  (3) Adapt the analysis for supervised learning problems to the current setting.
For two-layer neural network models, the function class is the Barron space \cite{EMaWu2019barron}.

At this point, only step (1) has been accomplished.

\begin{theorem}[A Posteriori Estimates \cite{Han2018convergence}]
  Under some assumptions, there exists a constant C, independent of h, d, and m, such that for sufficiently small h,
  \begin{equation}
  \begin{aligned}
    &\sup_{ t \in [0,T]} (\E|X_t - \hat{X}_t^{\pi}|^2 + \E|Y_t - \hat{Y}_t^{\pi}|^2) + \int_{0 }^{T} \E|Z_t - \hat{Z}_t^\pi|^2 \, \mathrm{d}t \\ 
    \le &C[h +  \E|g(X_{T}^\pi) - Y_T^{\pi}|^2],
    \end{aligned}
    \end{equation}
  where $\hat{X}_t^{\pi} = X_{t_i}^\pi$, $\hat{Y}_t^{\pi} = Y_{t_i}^\pi$, $\hat{Z}_t^{\pi} = Z_{t_i}^\pi$ for $t \in [t_i,t_{i+1})$.
\end{theorem}

\begin{theorem}[Upper Bound of Optimal Loss \cite{Han2018convergence}]
Under some assumptions, there exists a constant C, independent of h, d, and m, such that for sufficiently small h,
  \begin{equation}
  \begin{aligned}
      &\E|g(X_T^\pi) - Y_T^{\pi}|^2 \\
      \le~& C\Big\{h +  \E|Y_0 - \mu_0^\pi(\xi)|^2 + \sum_{i = 0}^{N-1}\E|\E[\tilde{Z}_{t_i}|X_{t_i}^\pi,Y_{t_i}^\pi] - \phi_i^\pi(X_{t_i}^\pi,Y_{t_i}^\pi)|^2h \Big\},
  \end{aligned}
  \end{equation}
where $\tilde{Z}_{t_i} =h^{-1}\E[\int_{t_i}^{t_{i+1}}Z_t\,\rmd t|\mathcal{F}_{t_i}]$. If $b$ and $\sigma$ are independent of $Y$, the term $\E[\tilde{Z}_{t_i}|X_{t_i}^\pi,Y_{t_i}^\pi]$ can be replaced with $\E[\tilde{Z}_{t_i}|X_{t_i}^\pi]$. 
\end{theorem}

\section{Control problems in high dimensions}
\label{sec:control}

One of the areas that high dimensional problems are often encountered is optimal control. In fact the term ``curse of dimensionality''
was first coined by Richard Bellman in the context of dynamic programming for control problems \cite{Bellman1957}.
Regarding CoD, there is an important difference between open- and closed-loop controls that we now explain.
f%
Consider the optimal control problem with a finite horizon $T$:

\begin{equation}
\label{eq:OCP}
\left \{
\begin{array}{cl}
\underset{ u}{\text{min}} & g ( x (T)) + \displaystyle \int_0^{T} L (t,  x(t), u(t)) dt , \\
\text{subject to} & \dot { x} (t) = f (t,  x,  u) , \\
        &  x (0) =  x_0 .
\end{array}
\right .
\end{equation}
Here $ x : [0, T] \to \mathcal X \subseteq \R^d$ is the state, $ u : [0, T] \to \mathcal U \subseteq  \R^m$ is the control, 
$g: \mathcal X \to \R$ is the terminal cost, and $L: [0, T \times \mathcal X \times \mathcal U \to \R$ is the running cost.
For fixed $\bx_0$, the problem above can be thought of as a two-point boundary value problem over the time interval $[0, T]$ and
the optimal control can be sought in the form: 
\begin{equation}
\label{eq: optimal openloop}
 u = u^* (t,  x^*(t))
\end{equation}
where $\bx^*$ denotes the optimal path.  We refer to \cite{Rao2009survey} for a review of the numerical algorithms for solving
this kind of two-point boundary value problems. In this case, CoD is not really an issue since the dimensionality of the independent
variable is just 1.
Controls of the form \eqref{eq: optimal openloop} is called an open-loop control.  In this case, the optimal control is only known along the optimal path.
Once the system deviates from the optimal path, one has to either recompute the optimal control or force
the system back to the optimal path.  
In many applications, one prefers a closed-loop control or feedback control
\begin{equation}
\label{eq: optimal feedback}
u = u^* (t, x),
\end{equation}
where the optimal control is known as every point in the state space. Closed-loop controls are functions of the state variable
and this is where the CoD problem arises.
To characterize open- and closed-loop controls, let 
\begin{equation}
\label{eq: Hamiltonian}
\tilde{ H} (t, x, \lambda, u) \coloneqq L (t, x, u) +  \lambda^T  f (t, x, u) 
\end{equation}
be the extended Hamiltonian associated with this control problem, and define
\begin{equation} 
\label{eq: optimal control as Hamiltonian minimizer}
 u^* (t, x; \lambda) = \underset{u \in \mathcal U}{\text{arg min}} \, H (t, x,  \lambda, u) .
\end{equation}
Here $\lambda$ is the co-state variable.  An important result is that the solution to the optimal control problem satisfies
Pontryagin's Minimum Principle:
\begin{equation}
\label{eq:BVP}
\begin{dcases*}
\dot{x} (t) = \frac{\del \tilde{H}}{\del \lambda} = f (t, x, u^* (t, x, \lambda)) , \\
\dot{\lambda} (t) = - \frac{\del \tilde{H}}{\del x} (t, x, \lambda, u^* (t, x, \lambda)) , \\
\dot v (t) = - L (t, x, u^* (t, x, \lambda)) ,
\end{dcases*}
\end{equation}
with the boundary conditions
\begin{equation}
\label{eq:BVP_BC}
x (0) = x_0 , \,
\lambda (T) = \nabla g(x (T)) , \,
v (T) = g (x (T)).
\end{equation}
Denote by $V$ the value function of the control problem:
\begin{equation}
\label{eq: value function}
V (t, x) \coloneqq \inf_{u \in \mathcal U} \left \{ g (y (T)) + \displaystyle \int_t^{T} L (\tau, y, u) d\tau \right \} ,
\end{equation}
subject to $\dot {y} (\tau) = f (\tau, y, u)$ and $y (t) = x$. 
Define the Hamiltonian:
\be
 H^* (t, x, \lambda) \coloneqq H (t, x, \lambda, u^*).
 \ee
 The HJB equation can be written as
\begin{equation}
\label{eq: HJB}
V_t (t, x) + H^* \left( t, x, V_{x} \right) = 0 
\end{equation}
with the terminal condition $V (T, x) = F (x)$.
The co-state and the closed-loop optimal control is given in terms of the value function by
\begin{equation}
\label{eq: costate as gradient}
\lambda (t) = \nabla_{\bx} V (t, x(t)),
\end{equation}
\begin{equation}
\label{eq: optimal feedback control using dVdx}
u^* (t, x) = \underset{u \in \mathcal U}{\text{arg min}} \, H \left( t, x, \nabla_{\bx} V, u \right) .
\end{equation}
To obtain an accurate approximation to the closed-loop control, we need to solve the control problem
for a large set of initial conditions, if not all.
The formulation \eqref{eq:OCP} is for a single initial condition.  To extend it to all initial conditions,
we consider instead the problem:
\begin{equation}
\label{eq: OCP1}
\underset{u}{\text{min}} \,  \E_{\bx_0 \sim \mu} \left( g (x (T)) + \displaystyle \int_0^{T} L (t, x(t), u(t, \bx(t))) dt \right)
\end{equation}
 subject to $\dot {x} (t) = f (t, x(t), u(t, (t)), x (0) = x_0 $.
 Here the optimization is over all possible policy functions $u$.
 One question that naturally arises is how we should choose the distribution $\mu$ for the initial condition.
 Clearly we are only interested in states whose value functions are not very big. Therefore one possible choice is
 the Gibbs distribution for the value function:
 \be
 \mu = \frac 1Z e^{-\beta V}
 \ee
 where $Z$ is a normalization factor. $\beta$ is a positive hyper-parameter. 

Unlike the stochastic case for which the training data is obtained on-the-fly, here one needs to address the issue of 
data generation explicitly. 
The following strategy was proposed in \cite{nakamura2019adaptive, kang2019algorithms}:
\bi
\item The two-point boundary value problem \eqref{eq:BVP}-\eqref{eq:BVP_BC} is solved to obtain the training data.
\item A neural network model is trained for the value function.
\ei
 In practice, \eqref{eq:BVP}-\eqref{eq:BVP_BC} is not an easy problem to solve, and
 it is important to  look for a small yet representative training dataset.
 The following ideas were proposed and tested in \cite{nakamura2019adaptive, kang2019algorithms}.
 
 The first is called ``warm start''.   The basic idea is to choose initializations for the iterative algorithms for \eqref{eq:BVP}-\eqref{eq:BVP_BC} to help guarantee convergence.  For example one can start with small values of $T$ in which case the convergence
 of the iterative algorithms is much less of an issue. One can use simple extrapolations of  these solutions on longer time intervals
 as initializations and obtain converged solutions on longer intervals.  This process can be continued.
 In addition,  once a reasonable approximation of the policy and value functions is obtained, one can use that to help initializing the
 two-point boundary value problem.
 
 The second is to explore adaptive sampling. %
 It has been explored in a similar context  \cite{zhang2020dp}.
 As for all adaptive algorithms, the key issue is an  error indicator: The larger the error, the more data are needed.
 \cite{zhang2020dp} uses the variance of the predictions from an ensemble of similar machine learning models as the error indicator,
 A  sophisticated error indicator that makes use of the variance of the gradient of the loss function was proposed in
 \cite{nakamura2019adaptive}.   Another idea is to simply use the magnitude of the gradient of the value function as the error indicator.

\section{Ritz,  Galerkin, and least squares}
The Ritz, Galerkin, and least squares formulations are among the most commonly
used frameworks for designing numerical algorithms for PDEs.
The Ritz formulation is based on a variational principle.
The Galerkin formulation is based on the weak formulation of a PDE that involves
 both the trial and test functions.
Least squares formulation is a very general approach for turning a PDE problem
into a variational problem by minimizing the squared residual of the PDE. 
It has the advantage of being general and straightforward to think about.
However, in classical numerical analysis, it is often the least preferred since the numerical
problem obtained this way tends to be worse conditioned than the ones using
Ritz or Galerkin formulation.
Designing machine learning-based algorithms using Ritz and least square formulations is rather
straightforward.  Since there is a variational principle behind both the Ritz and
least square formulations, one can simply replace the space of trial functions for
these variational principles by the hypothesis space in machine learning models.
Since machine learning is also a fundamentally optimization-based approach, the
integration of machine learning with variational methods for PDEs  is quite seamless.
Indeed these were among the first set of ideas that were proposed for machine learning-based
numerical algorithms for PDEs \cite{carleo2017solving,E2018deep,Sirignano2018dgm}.
For the same reason, designing machine learning-based algorithms using the Galerkin formulation is a different
matter, since Galerkin is not an optimization-based approach.  Rather it is based
on a weak formulation using test functions. 
The closest machine learning model to the Galerkin formulation is the Wasserstein GAN
(WGAN)
\cite{arjovsky2017towards,arjovsky2017wgan}: In WGAN, the discriminator plays the role of the test function;
the generator plays the role of the trial function. 
\subsection{The Deep Ritz method}
\begin{figure}
\centering
\includegraphics[width=0.48\textwidth,height=0.3\textwidth]{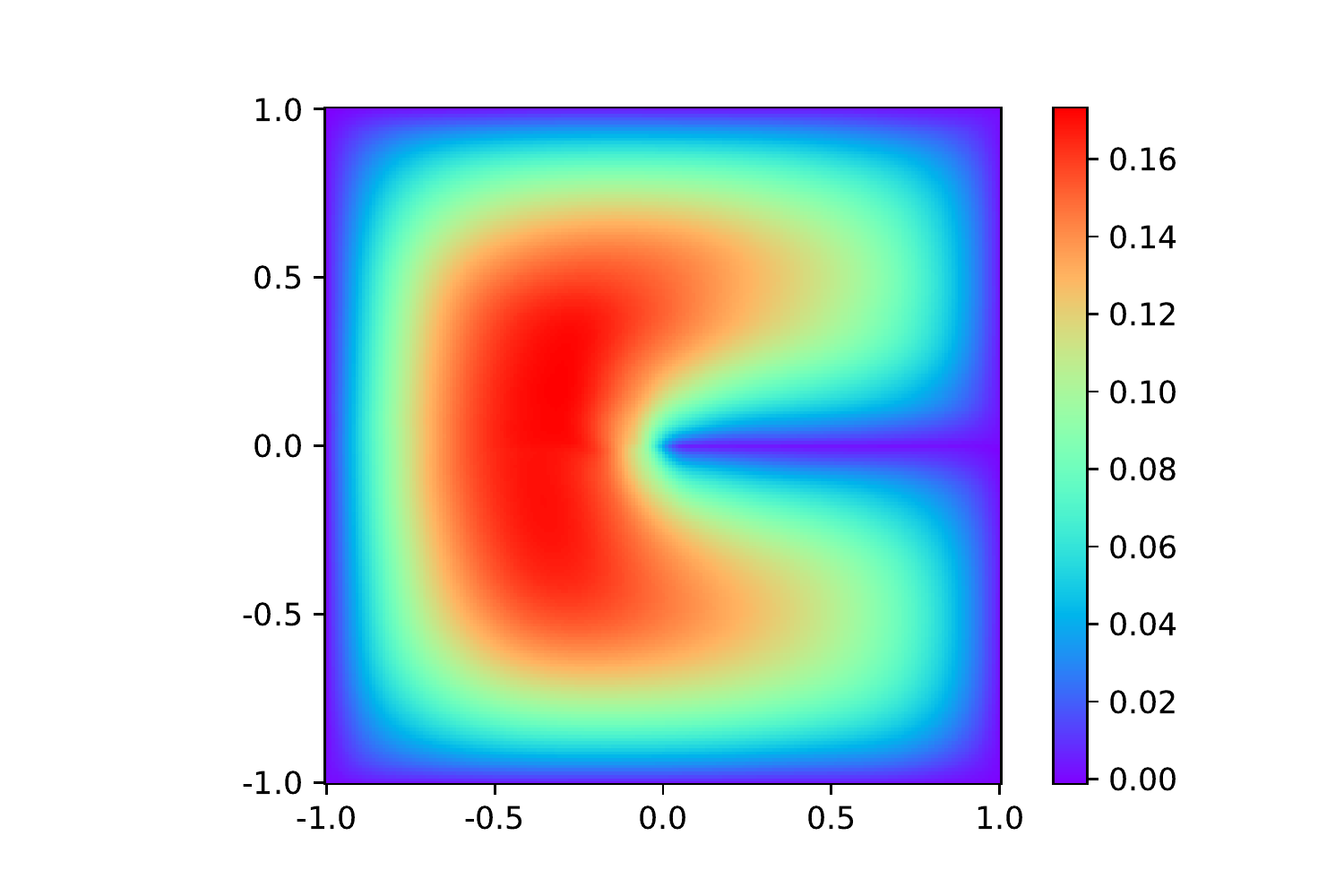}
\includegraphics[width=0.48\textwidth,height=0.3\textwidth]{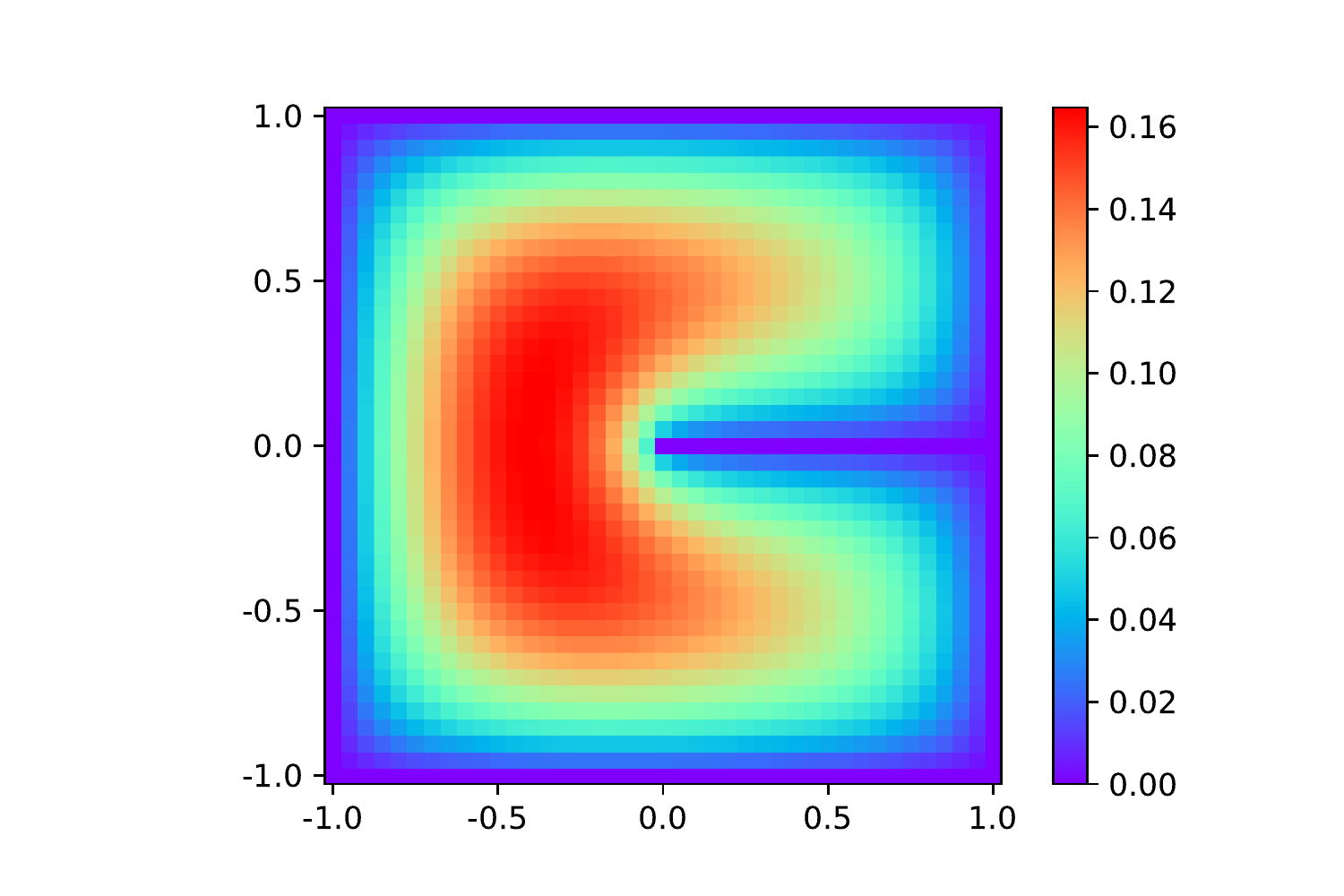}
\caption{Solutions computed by two different methods. On the left is Deep Ritz with
811 parameters. On the right is the solution of the finite difference method on a uniform
grid with 1681 parameters.
Reprinted from \cite{E2018deep}.
}
\label{fig:ritz}
\end{figure}

The Deep Ritz method was proposed in \cite{E2018deep}.
Consider the variational problem \cite{Evans2010partial}
\begin{equation}
\min_{u \in H} I(u)
\end{equation}
where
\begin{equation}
I(u) = \int_{\Omega} \left( \frac 12  |\nabla u(x)|^2 -  f(x) u(x) \right) dx
\label{I(u)}
\end{equation}
and $H$ is the set of admissible functions (also called
trial function, here represented by $u$),
$f$ is a given function, representing
external forcing to the system under consideration.
It is understood that boundary conditions are incorporated into the definition of $H$.
The Deep Ritz method consists of the following components:
\begin{enumerate}
\item Deep neural network-based representation of the trial function.
\item A numerical quadrature rule for the functional.
\item An algorithm for solving the final optimization problem.
\end{enumerate}
Each component is relatively straightforward. One can take the usual neural network
models to represent the trial function. 
In high dimensions one needs an effective 
 Monte Carlo algorithm to discretize the integral in \eqref{I(u)}.
The interplay between the discretization of the integral and the discretization of the trial
function using neural network models is an interesting issue that requires further attention.
Finally, SGD can be used naturally, similar to the situation in Deep BSDE: The integral
in the functional \eqref{I(u)} plays the role of the expectation in Deep BSDE.
One notable issue is the choice of the activation function. ReLU activation does not perform well
due to the discontinuity in its derivative.
It has been observed that the activation function $\sigma_3 (z) = \max (z, 0)$ performs
 much better than ReLU.  More careful study is needed on this issue also.
One feature of the Deep Ritz method that potentially makes it interesting even
for low dimensional problems is that it is mesh-free and naturally adaptive.
To examine this we consider the well-known crack problem: Computing the displacement 
around a crack.  To this end, we consider the Poisson equation:
\begin{equation}
\begin{aligned}
-&\Delta u(x)=1,\quad &x\in \Omega\\
&u(x)=0,\quad &x\in \partial \Omega
\end{aligned}
\end{equation}
where $\Omega=(-1,1)\times (-1,1) \backslash [0,1)\times \{0\}$.
The solution to this problem suffers from the well-known ``corner singularity"
caused by the nature of the domain \cite{strang1973analysis}.
A simple asymptotic analysis shows that at the origin,  the solution
behaves as
$u(x)=u(r,\theta) \sim r^{\frac{1}{2}}\sin\frac{\theta}{2}$
\cite{strang1973analysis}.
Models of this type have been extensively used
to help developing and testing adaptive finite element methods.
Here the essential boundary condition causes some problems.
The simplest idea is to just  use a penalty method and consider
the modified functional
\begin{equation}
I(u) = \int_\Omega \left( \frac{1}{2} |\nabla_x u(x)|^2 - f(x) u(x) \right) dx
+\beta \int_{\partial \Omega} u(x)^2 ds.
\end{equation}
An acceptable choice is $\beta = 500$.
The results from the Deep Ritz method 
with $811$ parameters in the neural network model
and the finite difference method with $\Delta x_1=\Delta x_2=0.1$ ($1,681$
degrees of freedom) are shown in Figure \ref{fig:ritz}.
More quantitative comparisons can be found in \cite{E2018deep}.
Of course adaptive numerical methods are very well developed for solving
problems with corner singularities and even more general singular problems.
Nevertheless, this example shows that Deep Ritz is potentially a naturally
adaptive algorithm.
There are also a number of problems that need to be addressed in future work:
\begin{enumerate}
\item  The variational problem that results from Deep Ritz is usually not
convex even when the original problem is.
\item At the present time, there are no consistent conclusions about the
convergence rate.
\item The treatment of the essential boundary condition
is not as simple as the traditional methods.
\end{enumerate}
Some analysis of the  Deep Ritz method has been carried out in \cite{muller2019deep}.
\subsection{The least square formulation}
The least square approach was used in \cite{carleo2017solving} for solving the
dynamic Schr\"odinger equation and was subsequently developed more systematically
in \cite{Sirignano2018dgm} (although \cite{Sirignano2018dgm} referred to it as Galerkin method).
The basic idea is very simple: Solving the PDE
\be
Lu = f
\label{PDE-1}
\ee
over a domain $\Omega$ in $\R^d$ can be formulated equivalently as solving the
variational problem for the functional
\be
J(u) = \int_{\Omega} \| Lu - f \|^2 \mu(d \bx)
\ee
where $\mu$ is a suitably chosen probability distribution on $\Omega$.
$\mu$ should be non-degenerate and readily sampled.
With this,  the least square formulation looks very similar to the Ritz formulation with $J$ replacing
the functional $I$.
\subsection{The Galerkin formulation}
The starting point of the Galerkin approximation is the weak form of 
\eqref{PDE-1}:
\be
a(u, \phi)  =  (Lu, \phi) = (f, \phi), u \in H_1, \phi \in H_2
\label{Galerkin-1}
\ee
where $H_1$ and $H_2$ are the trial and test function spaces respectively,
$\phi$ is an arbitrary test function in $H_2$, 
$(\cdot, \cdot) $ is the standard $L^2$ inner product for functions. 
Usually some integration by parts is applied.  For example, if $L = - \Delta$, then
except for boundary terms, one has
\be
a(u, \phi) = (\nabla u, \nabla \phi) 
\ee
Therefore this formulation only involves first order derivatives.
The most important feature of the Galerkin formulation is that involves
the test function.
In this spirit, the Wasserstein GAN can also be regarded as an example of the
Galerkin formulation. Given a set of data $\{\bx_j, j=1, 2, \cdots ,n\}$ in
$\R^d$ and a reference probability distribution $\nu^*$ on $\R^{d'}$,
we look for the mapping $G$ (the generator) from $\R^{d'}$ to $\R^d$,
such that \cite{arjovsky2017wgan}
\be
\int_{\R^{d'}} \phi(G(\bz)) \nu^*(d \bz) = \frac 1n \sum_{j=1}^n \phi(\bx_j)
\label{WGAN}
\ee
for all Lipschitz functions $\phi$. The test function $\phi$ is called the discriminator
in this context.
Like GAN, the most obvious reformulation of \eqref{Galerkin-1} is a min-max problem:
\be
\min_{u \in H_1} \max_{\|\phi\|_{H_2} \le 1} (a(u, \phi) - (f, \phi))^2.
\label{Galerkin-2}
\ee
Unfortunately this formulation is not easy to work with.
The problems encountered are similar to those in WGAN. 
Despite this some very encouraging progress has been made and we
 refer to \cite{Zang2020weak} for the details.

\section{Multilevel Picard approximation methods for nonlinear PDEs}
\label{sec:MLP}

In the articles 
E et al.~\cite{E2016multilevel}
and 
Hutzenthaler et al.~\cite{Hutzenthaleretal2018arXiv}
so-called
fully history recursive multilevel Picard approximation methods 
have been introduced and analyzed (in the following we abbreviate 
\emph{fully history recursive multilevel Picard} by MLP). 
The error analysis in the original article 
Hutzenthaler et al.~\cite{Hutzenthaleretal2018arXiv}
is restricted to semilinear heat PDEs with Lipschitz nonlinearities. 
By now in the scientific literature there are, however, a series of 
further articles on such MLP approximation methods 
(see \cite{hutzenthaler2019arxiv1903,Becketal2019MLP_nonlip_arXiv,giles2019generalised,beck2020arxiv2003,becker2020arxiv2005,hutzenthaler2019arxiv1912,hutzenthaler2020multilevel,
E2019multilevel,hutzenthaler2020lipschitz}) 
which analyze, extend, or generalize the 
MLP approximation methods proposed in \cite{E2016multilevel,Hutzenthaleretal2018arXiv} 
to larger classes of PDE problems such as 
semilinear Black-Scholes PDEs (see \cite{hutzenthaler2019arxiv1903,becker2020arxiv2005}), 
semilinear heat PDEs with gradient dependent nonlinearities (see \cite{hutzenthaler2019arxiv1912,hutzenthaler2020multilevel}), 
semilinear elliptic PDE problems (see \cite{beck2020arxiv2003}), 
semilinear heat PDEs with non-Lipschitz continuous nonlinearities (see \cite{Becketal2019MLP_nonlip_arXiv,becker2020arxiv2005}), 
and 
semilinear second-order PDEs with varying coefficient functions (see~\cite{hutzenthaler2019arxiv1903,hutzenthaler2020lipschitz}).

In the remainder of this section we sketch the main ideas of MLP approximation methods 
and to keep the presentations as easy as possible we restrict ourselves in the following 
presentations to semilinear heat PDEs with Lipschitz continuous nonlinearities 
with bounded initial values. 
The next result, Theorem~\ref{thm:MLP_analysis} below, provides a complexity analysis 
for MLP approximations in the case of semilinear heat PDEs with Lipschitz continuous nonlinearities. 
Theorem~\ref{thm:MLP_analysis} is strongly based 
on Hutzenthaler et al.~\cite[Theorem~1.1]{Hutzenthaleretal2018arXiv} and 
Beck et al.~\cite[Theorem~1.1]{Becketal2019MLP_nonlip_arXiv}. 
\begin{theorem}
\label{thm:MLP_analysis}
Let
$ T \in (0,\infty) $,
$ \Theta = \cup_{ n \in \N } \Z^n $, 
let $ f \colon \R \to \R $ be Lipschitz continuous, 
for every $ d \in \N $ let 
$ u_d \in C^{ 1, 2 }( [0,T] \times \R^d, \R ) $
be at most polynomially growing, 
assume for every 
$ d \in \N $, $ t \in [0,T] $,
$ x \in \R^d $
that 
\begin{equation}
\label{eq:MLP_allen_cahn_pde}
  (\tfrac{\partial}{\partial t}u_d)(t,x) 
  = 
  (\Delta_x u_d)(t,x) 
  + 
  f(u_d(t,x)),
\end{equation} 
let 
  $(\Omega,\mathcal{F},\P)$ 
be a probability space, 
let
  $\mathcal{R}^{\theta}\colon \Omega \to [0,1]$, $\theta\in\Theta$, 
be independent $\mathcal{U}_{[0,1]}$-distributed random variables, 
let 
  $W^{d,\theta}\colon [0,T]\times\Omega\to\R^d$, 
  $d\in\N$,
  $\theta\in\Theta$,
be independent standard Brownian motions, 
assume that 
  $(\mathcal{R}^{\theta})_{\theta\in\Theta}$ 
and
  $(W^{d,\theta})_{(d,\theta)\in\N\times\Theta}$ 
are independent, 
for every 
  $d\in \N$,
  $s\in [0,T]$, 
  $t\in [s,T]$,
  $x\in \R^d$, 
  $\theta\in \Theta$ 
let 
  $X^{d,\theta}_{s,t,x}\colon\Omega\to\R^d$
satisfy 
  $
  X^{d,\theta}_{s,t,x}
  = 
  x + \sqrt{2}(W^{d,\theta}_t - W^{d,\theta}_s)
  $, 
let 
  $U^{d,\theta}_{n,M}\colon [0,T]\times\R^d\times\Omega\to\R$, 
  $d,n,M\in\N_0$,
  $\theta\in\Theta$, 
satisfy for every 
  $d,M\in\N$,
  $ n \in \N_0 $, 
  $\theta\in\Theta$,
  $t\in [0,T]$, 
  $x\in \R^d$ 
that
  \begin{align}\label{eq:mlp}
  & 
  U^{d,\theta}_{n,M}(t,x) 
  = 
  \sum_{k=1}^{n-1} \frac{t}{M^{n-k}} 
  \Bigg[ 
  \sum_{m=1}^{M^{n-k}} 
  \bigg(
  f\Big( 
  U^{d,(\theta,k,m)}_{k,M}\big( t \mathcal{R}^{(\theta,k,m)} , 
  X^{d,(\theta,k,m)}_{ t \mathcal{R}^{(\theta,k,m)},t,x}
  \big)
  \Big)
  \\ 
  & 
  -
  f\Big( 
  U^{d,(\theta,-k,m)}_{k-1,M}\big( t \mathcal{R}^{(\theta,k,m)} , 
  X^{d,(\theta,k,m)}_{ t \mathcal{R}^{(\theta,k,m)},t,x}\big)
  \Big)
  \bigg)
  \Bigg]
  + 
  \frac{ \mathbbm{1}_{ \N }( n ) }{M^n}\!\left[ 
  \sum_{m=1}^{M^n} 
  \left( 
  u_d(0,X^{d,(\theta,0,-m)}_{0,t,x}) 
  + 
  t \, f(0)
  \right)
  \right], 
  \nonumber
  \end{align} 
and for every 
$ d, M \in \N $, 
$ n \in \N_0 $
let 
$ \mathfrak{C}_{ d, n, M } \in \N_0 $
be the number of function evaluations of $ f $ and $ u_d(0,\cdot) $ and 
the number of realizations of scalar random variables 
which are used to compute one realization of 
$ U^{d,0}_{n,M}(T,0) \colon \Omega \to \R $ 
(cf.~\cite[Corollary 4.4]{hutzenthaler2020lipschitz} for a precise definition). 
Then there exist
$ \mathfrak{N} \colon (0,1] \to \N $ 
and 
$ c \in \R $ 
such that for all 
$ d \in \N $, $ \varepsilon \in (0,1] $ 
it holds that 
$
  \mathfrak{C}_{d,\mathfrak{N}_{\varepsilon},\mathfrak{N}_{\varepsilon}} 
  \leq c d^c \varepsilon^{ - 3 }
$
and 
$ 
  \big(  
    \E\big[ 
      |
        U^{d,0}_{\mathfrak{N}_{\varepsilon},\mathfrak{N}_{\varepsilon}}(T,0)
        - 
        u_d(T,0)
      |^2
    \big]
  \big)^{\!\nicefrac12} 
  \leq \varepsilon
$. 
\end{theorem} 
In the following we add some comments on the statement in Theorem~\ref{thm:MLP_analysis} above 
and we thereby also provide explanations for some of the mathematical objects which appear in Theorem~\ref{thm:MLP_analysis}.

Theorem~\ref{thm:MLP_analysis} provides a complexity analysis for MLP approximations 
in the case of semilinear heat PDEs with Lipschitz continuous nonlinearities. 
In \eqref{eq:mlp} in Theorem~\ref{thm:MLP_analysis} the employed MLP approximations are specified. 
The MLP approximations in \eqref{eq:mlp} aim to approximate the solutions of the PDEs in \eqref{eq:MLP_allen_cahn_pde}. 
The strictly positive real number $ T \in (0,\infty) $ in Theorem~\ref{thm:MLP_analysis} describes the time horizon of the PDEs in \eqref{eq:MLP_allen_cahn_pde}. The function $ f \colon \R \to \R $ in Theorem~\ref{thm:MLP_analysis} describes \emph{the nonlinearity of the PDEs} in \eqref{eq:MLP_allen_cahn_pde}. For simplicity we restrict ourselves in 
Theorem~\ref{thm:MLP_analysis} in this article to Lipschitz continuous 
nonlinearities which do only depend on the solution of the PDE but not 
on the time variable $ t \in [0,T] $, not on the space variable $ x \in \R^d $, and also not on the derivatives 
of the PDE solution. In the more general MLP analyses in the scientific literature 
(cf., 
e.g., \cite{
Hutzenthaleretal2018arXiv,
hutzenthaler2019arxiv1903,
Becketal2019MLP_nonlip_arXiv,
giles2019generalised,
beck2020arxiv2003,
becker2020arxiv2005,
hutzenthaler2019arxiv1912,
hutzenthaler2020lipschitz}) the nonlinearity of the PDE is allowed to depend on the time variable $ t \in [0,T] $, 
on the space variable $ x \in \R^d $, 
on the PDE solution, 
and also on the derivatives of the PDE solution (see~\cite{hutzenthaler2019arxiv1912}), 
and the nonlinearity of the PDE may also be not Lipschitz continuous (see~\cite{Becketal2019MLP_nonlip_arXiv,becker2020arxiv2005}).

The functions $ u_d \colon [0,T] \times \R^d \to \R $, $ d \in \N $, 
in Theorem~\ref{thm:MLP_analysis} 
describe the exact solutions of the PDEs in \eqref{eq:MLP_allen_cahn_pde}. 
The linear differential operator on the 
right hand side of the PDE in 
\eqref{eq:MLP_allen_cahn_pde} 
is just the Laplacian and 
Theorem~\ref{thm:MLP_analysis} thus only applies to 
semilinear heat PDEs of the form \eqref{eq:MLP_allen_cahn_pde} 
but the MLP analyses in the scientific literature 
also apply to PDEs with 
\emph{much more general second-order differential operators} 
(cf., e.g.,~\cite{hutzenthaler2020lipschitz,
hutzenthaler2019arxiv1903}).

The approximation algorithm in \eqref{eq:mlp} is a Monte-Carlo algorithm in the sense 
that it employs Monte-Carlo averages based on many independent identically distributed (i.i.d.) 
random variables. In the case of plain-vanilla standard Monte Carlo algorithms 
for linear PDEs the employed i.i.d.\ random variables are often indexed 
through the set of all natural numbers $ \N = \{ 1, 2, 3, \dots \} $ 
where we have one random variable for each natural number $ n \in \N $. 
The approximation algorithm in \eqref{eq:mlp} is somehow a \emph{nonlinear Monte-Carlo algorithm} 
and in the case of such a nonlinear Monte-Carlo algorithm the situation is getting more 
complicated and, roughly speaking, we need more i.i.d.\ random variables and 
therefore, roughly speaking, also a larger index set\footnote{We remark that the set of all natural numbers $ \N $ 
is a proper subset of the set $ \Theta = \cup_{ n \in \N } \Z^n $ in the sense that $ \N \subsetneqq \Theta $ but 
the set of all natural numbers $ \N $ 
and the set $ \Theta $ have, of course, the same cardinality}. More precisely, 
in the case of the MLP approximation algorithm in \eqref{eq:mlp} we employ the 
set $ \Theta = \cup_{ n \in \N } \Z^n $ in Theorem~\ref{thm:MLP_analysis} 
as the index set to introduce sufficiently many i.i.d.\ random variables. 
In particular, in Theorem~\ref{thm:MLP_analysis} we use 
the family $ \mathcal{R}^{ \theta } \colon \Omega \to [0,1] $,
$ \theta \in \Theta $, of independent 
on $ [0,1] $ continuous uniformly distributed random variables 
and the family 
$ W^{ d, \theta } \colon [0,T] \times \Omega \to \R^d $, 
$ d \in \N $,
$ \theta \in \Theta $, 
of independent standard Brownian motions 
as random input sources for the MLP approximation algorithm 
in \eqref{eq:mlp}.

The natural numbers 
$ \mathfrak{C}_{ d, n, M } \in \N_0 $, 
$ d, M \in \N $, 
$ n \in \N_0 $,
in Theorem~\ref{thm:MLP_analysis} above 
aim to measure to 
\emph{computational cost} 
for the MLP approximation algorithm in \eqref{eq:mlp}. 
Theorem~\ref{thm:MLP_analysis} shows that 
there exists a function 
$ \mathfrak{N} \colon (0,1] \to \N $ 
and a real number $ c \in \R $ 
such that for all 
$ d \in \N $, $ \varepsilon \in (0,1] $ 
we have that 
the $ L^2 $-approximation error 
$
  \big(  
    \E\big[ 
      |
        U^{d,0}_{\mathfrak{N}_{\varepsilon},\mathfrak{N}_{\varepsilon}}(T,0)
        - 
        u_d(T,0)
      |^2
    \big]
  \big)^{\!\nicefrac12} 
$
between the MLP approximation 
$
  U^{d,0}_{\mathfrak{N}_{\varepsilon},\mathfrak{N}_{\varepsilon}}(T,0)
$
and the exact solution 
$ u_d $
of the PDE 
at the space-time point 
$ t = T $, $ x = 0 $
is smaller or equal than the 
prescribed approximation 
accuracy $ \varepsilon $ 
with the computational cost 
$
  \mathfrak{C}_{d,\mathfrak{N}_{\varepsilon},\mathfrak{N}_{\varepsilon}} 
$
for the MLP approximations 
being smaller or equal than 
$
  c d^c \varepsilon^{ - 3 }
$. 
The computational cost 
for the MLP approximation algorithm
thus grows at most polynomially 
in the PDE dimension $ d \in \N $ 
and at most cubically 
in the reciprocal $ \varepsilon^{ - 1 } $ 
of the prescribed approximation accuracy 
$ \varepsilon \in (0,1] $.

The more general MLP approximation
\cite{Becketal2019MLP_nonlip_arXiv,
Hutzenthaleretal2018arXiv,
hutzenthaler2019arxiv1903,
giles2019generalised,
hutzenthaler2019arxiv1912,
hutzenthaler2020lipschitz}
in the scientific literature 
improve this statement in several ways. 
First, the main approximation results in the above named reference list 
allow the numerical approximation of the PDE solution not 
necessarily at the space point $ x = 0 $ 
but at much more general space points. 
Second, the main approximation results in the above named reference list
provide explicit errors constants and explicit exponents in dependence on the constants in 
the assumptions for the involved functions. 
For instance, if the initial conditions 
of the PDEs under consideration are bounded functions, then 
Hutzenthaler et~al.~\cite[Theorem~1.1]{Hutzenthaleretal2018arXiv} and 
Beck et~al.~\cite[Theorem~1.1]{Becketal2019MLP_nonlip_arXiv} 
even prove that the computation cost of the employed MLP approximation scheme 
\emph{grows at most linearly in the PDE dimension}. 
Finally, most of the MLP approximation results in the scientific literature also prove 
that the computational cost of the considered MLP approximation scheme 
grows up to an arbitrarily small real number 
at most quadratically (instead of cubically as in Theorem~\ref{thm:MLP_analysis} above) 
in the reciprocal $ \varepsilon^{ - 1 } $
of the prescribed approximation accuracy 
$ \varepsilon \in (0,1] $.

It should also be noted that MLP approximation 
schemes not only overcome the curse of dimensionality 
in the numerical approximation of parabolic PDEs but also in the case of \emph{elliptic PDEs} with Lipschitz 
continuous nonlinearities (see Beck et al.~\cite{beck2020arxiv2003}). 
Encouraging numerical simulations for MLP approximation schemes 
in the case of semilinear Black-Scholes PDEs, systems of semilinear PDEs, Allen-Cahn PDEs, 
and sine-Gordon type PDEs can be found in Becker et al.~\cite{becker2020arxiv2005} (see also E et al.~\cite{E2019multilevel}).

In this article we do not provide a detailed proof for Theorem~\ref{thm:MLP_analysis} 
but instead we refer, e.g., Hutzenthaler et al.~\cite{Hutzenthaleretal2018arXiv,hutzenthaler2020lipschitz}
for a detailed proof of Theorem~\ref{thm:MLP_analysis}. 
In addition, in the following we also briefly outline 
some of \emph{the main ideas of the proof} of Theorem~\ref{thm:MLP_analysis}. 
The derivation and thereby also the mathematical analysis of the MLP approximation schemes 
is, roughly speaking, based on the following three steps. 
\begin{enumerate}[(I)]
\item 
First, we reformulate the PDE under consideration (or, more generally, the computational problem under consideration) 
as a suitable stochastic fixed point equation with the unique fixed point solution of 
the stochastic fixed point equation being the unique solution of the PDE under consideration. 
\item 
Second, we approximate the unique fixed point of the stochastic fixed point equation by means 
of fixed point iterations according to the Banach fixed point theorem (which are referred to as 
Picard iterations in the context of temporal integral fixed point equations). 
\item 
Third, we recursively approximate the resulting fixed point iterations by means of 
suitable multilevel Monte-Carlo approximations resulting with the resulting Monte-Carlo approximations 
being full history recursive. 
\end{enumerate}
A key idea in the above derivation of the MLP approximation scheme is that the fixed point iterations often converge 
exceedingly quick, that is, with factorial convergence speed to the unique fixed point of the stochastic fixed 
point equation under consideration while efficient multilevel Monte Carlo approximations assure that the 
computation cost of the considered MLP approximation scheme grows 
not significantly larger than factorially. 
These facts made it possible to prove that MLP approximation schemes overcome the curse of dimensionality 
in the numerical approximation of a large class of semilinear PDEs.

Despite the great performance of deep learning-based approximation schemes in various numerical simulations, 
until today, MLP approximation schemes 
are, to the best of our knowledge, 
the only approximation schemes 
in the scientific literature 
for which 
\emph{it has been proven that they 
do indeed overcome the curse of dimensionality}  
in the numerical approximation of semilinear PDEs with general time horizons.

\section{Mathematical results for neural network approximations for PDEs}

Until today, there is no complete rigorous mathematical analysis 
which proves (or disproves) the conjecture that there exists 
a deep learning-based approximation method which overcomes the curse 
of dimensionality in the numerical approximation of PDEs. 
However, there are now a few mathematical results in the scientific literature 
(see, e.g.,~\cite{berner2018analysis,elbrachter2018dnn,gonon2019uniform,grohs2018proof,grohs2019deep,hutzenthaler2020proof,jentzen2018proofarxiv1809,kutyniok2019theoretical,reisinger2019rectified,HornungJentzenSalimova2020}) which prove that deep neural networks have the capacity to approximate 
solutions of PDEs without the curse of dimensionality.

In particular, in the article Grohs et al.~\cite{grohs2018proof} it has been proved 
that there exist neural networks which approximate solutions of linear 
Black-Scholes PDEs with the number of parameters of the neural networks 
growing at most polynomially in both 
the reciprocal $ \nicefrac{ 1 }{ \varepsilon } $ 
of the prescribed approximation accuracy $ \varepsilon \in (0,\infty) $
and the PDE dimension $ d \in \N $. 
The articles~\cite{berner2018analysis,elbrachter2018dnn,gonon2019uniform,grohs2019deep,jentzen2018proofarxiv1809,kutyniok2019theoretical,reisinger2019rectified,HornungJentzenSalimova2020}, 
in particular, extend 
the results in the article Grohs et al.~\cite{grohs2018proof}
to more general linear PDEs and the article 
Hutzenthaler et al.~\cite{hutzenthaler2020proof} 
extends the results in the article Grohs et al.~\cite{grohs2018proof}
to nonlinear heat PDEs with Lipschitz continuous nonlinearities. 
To better explain the results in the article 
Hutzenthaler et al.~\cite{hutzenthaler2020proof}, 
we now present in the following result, 
Theorem~\ref{thm:DNN_nonlinear_PDEs} below, 
a special case of 
Hutzenthaler et 
al.~\cite[Theorem~1.1]{hutzenthaler2020proof}. 
\begin{theorem}
\label{thm:DNN_nonlinear_PDEs}
Let $ \rho \colon (\bigcup_{d\in\N}\R^d) \to (\bigcup_{d\in\N}\R^d)$ 
satisfy for all $d\in\N $, $ x = (x_1,\dots,x_d) \in \R^d $ 
that 
$
  \rho(x) = ( \max\{x_1,0\}, \allowbreak \dots, \allowbreak\max\{x_d,0\} ) 
$, 
let $
  {\bf N} = 
  \bigcup_{ L \in \N }
  \bigcup_{ l_0, l_1, \dots, l_L \in \N }(
    \bigtimes_{ k = 1 }^L ( \R^{ l_k \times l_{k-1} } \times \R^{l_k} )
  ) 
$, 
let 
$ \mathcal{R} \colon {\bf N} \to (\bigcup_{k,l\in\N}C(\R^k,\R^l)) $ 
and 
$
  \mathcal{P} \colon {\bf N} \to \N 
$ 
satisfy for all $L\in\N,$ $l_0,l_1,\dots,l_L\in\N,$ $\Phi = ((W_1,B_1),\allowbreak(W_2,B_2),\allowbreak\dots,\allowbreak(W_L,B_L)) \in (\bigtimes_{k=1}^L\allowbreak (\R^{l_k\times l_{k-1}}\times \R^{l_k})),$ $x_0\in \R^{l_0},$ $x_1\in\R^{l_1},$ $\dots,$ $x_L\in\R^{l_L}$ with $\forall\, k\in \{1,2,\dots,L-1\}\colon x_k = \rho( W_k x_{ k - 1 } + B_k ) $ that
$
  \mathcal{R}(\Phi) \in C(\R^{l_0},\R^{l_{L}})
$, 
$
  (\mathcal{R}(\Phi))(x_0) = W_Lx_{L-1}+B_L
$, 
and
$
  \mathcal{P}(\Phi) = \sum_{k=1}^L l_k(l_{k-1}+1)
$, 
let $ T, \kappa \in (0,\infty) $, 
$ (\mathfrak{g}_{d,\varepsilon} )_{ (d,\varepsilon) \in \N \times (0,1] } \subseteq {\bf N} $, 
let $f\colon\R \to \R$ be Lipschitz continuous, 
let $u_d \in C^{1,2}([0,T]\times \R^d,\R),$ $d\in\N$, 
and assume for all 
$ d \in \N $, $ x = (x_1, \dots, x_d) \in \R^d $, 
$ \varepsilon \in (0,1] $, $ t \in [0,T] $ 
that 
$
  \mathcal{R}(\mathfrak{g}_{d,\varepsilon}) \in C(\R^d,\R) 
$, 
$ 
  \varepsilon | u_d(t,x) | 
  + \allowbreak 
  | u_d(0,x) \allowbreak - ( \mathcal{R}(\mathfrak{g}_{d,\varepsilon}) )(x) | 
  \le \varepsilon \kappa d^\kappa (1 + \sum_{ i = 1 }^d | x_i |^\kappa ) 
$, 
$
  \mathcal{P}(\mathfrak{g}_{d,\varepsilon}) \le \kappa d^\kappa \varepsilon^{-\kappa} 
$, 
and 
\begin{equation}\label{eq:1a}
  (\tfrac{\partial}{\partial t} u_d)(t,x) = (\Delta_x u_d)(t,x) + f(u_d(t,x)) .
\end{equation}
Then there exist 
$ (\mathfrak{u}_{d,\varepsilon})_{(d,\varepsilon)\in \N\times(0,1] }\subseteq {\bf N} $
and $ c \in \R $ 
such that for all $ d \in \N $, $ \varepsilon \in (0,1] $ 
it holds that 
$ \mathcal{R}( \mathfrak{u}_{d,\varepsilon} ) \in C( \R^d, \R ) $, 
$ \mathcal{P}( \mathfrak{u}_{d,\varepsilon} ) \le c d^c \varepsilon^{-c} $, 
and
\begin{equation}\label{eq:2a}
  \left[ 
    \int_{ [0,1]^d } 
      \left| 
        u_d(T,x) - (\mathcal{R}(\mathfrak{u}_{d,\varepsilon}))(x) 
      \right|^2 
    dx
  \right]^{\nicefrac{1}{2}} \le \varepsilon.
\end{equation}
\end{theorem}
Theorem~\ref{thm:DNN_nonlinear_PDEs} is an immediate 
consequence of 
Hutzenthaler et al.~\cite[Theorem~1.1]{hutzenthaler2020proof}. 
In the following we add some comments 
on the mathematical objects 
appearing in Theorem~\ref{thm:DNN_nonlinear_PDEs} above and, thereby, we 
also add some explanatory comments on the statement of Theorem~\ref{thm:DNN_nonlinear_PDEs}.

Theorem~\ref{thm:DNN_nonlinear_PDEs} is a DNN approximation result 
with the activation functions in the DNNs being multidimensional rectifier functions 
described by the function 
$ \rho \colon (\bigcup_{d\in\N}\R^d) \to (\bigcup_{d\in\N}\R^d) $ 
in Theorem~\ref{thm:DNN_nonlinear_PDEs} above. 
The set $ {\bf N} $ in Theorem~\ref{thm:DNN_nonlinear_PDEs} above 
represents \emph{the set of all neural networks}. 
The function 
$ 
  \mathcal{R} \colon {\bf N} \to 
  ( \bigcup_{ k, l \in \N } C( \R^k, \R^l ) ) 
$ 
maps neural networks to their realization function in the sense that for every 
$ \Phi \in {\bf N} $ 
it holds that 
$
  \mathcal{R}( \Phi )
  \in 
  ( \bigcup_{ k, l \in \N } C( \R^k, \R^l ) ) 
$
is the realization function associated with the neural network $ \Phi $. 
The function 
$
  \mathcal{P} \colon {\bf N} \to \N 
$ 
counts the number of parameters of the neural networks 
in the sense that for every 
$ \Phi \in {\bf N} $ 
it holds that 
$
  \mathcal{P}( \Phi ) \in \N
$
represents the number of real numbers which are used 
to uniquely describe the neural network $ \Phi $.

Theorem~\ref{thm:DNN_nonlinear_PDEs} demonstrates that the solutions 
of the PDEs in \eqref{eq:1a} can be approximated by DNNs 
without the curse of dimensionality. 
The real number $ T \in (0,\infty) $ in Theorem~\ref{thm:DNN_nonlinear_PDEs} describes 
the time horizon of the PDEs in \eqref{eq:1a} above. 
The function $ f \colon \R \to \R $ in Theorem~\ref{thm:DNN_nonlinear_PDEs} 
describes the nonlinearity 
in the PDEs in \eqref{eq:1a}. 
It is assumed to be Lipschitz continuous in the sense that there exists 
$ L \in \R $ such that for all $ x, y \in \R $ it holds that 
$ | f(x) - f(y) | \leq L | x - y | $.

The real number $ \kappa \in (0,\infty) $ 
in Theorem~\ref{thm:DNN_nonlinear_PDEs}  
is used to formulate 
\emph{the regularity and 
the approximation assumptions} 
which we impose in Theorem~\ref{thm:DNN_nonlinear_PDEs}. 
In particular, we assume in Theorem~\ref{thm:DNN_nonlinear_PDEs} 
that the initial value functions of the PDEs in \eqref{eq:1a} 
can be approximated by DNNs without the curse of dimensionality. 
In Theorem~\ref{thm:DNN_nonlinear_PDEs} this approximation assumption 
is formulated by means of the family 
$ 
  (\mathfrak{g}_{d,\varepsilon} )_{ (d,\varepsilon) \in \N \times (0,1] } 
  \subseteq {\bf N} 
$
of neural networks. 
More formally, 
in Theorem~\ref{thm:DNN_nonlinear_PDEs} 
we assume that there exist neural networks   
$ \mathfrak{g}_{d,\varepsilon} \in {\bf N} $,
$ 
  d \in \N 
$,
$
  \varepsilon \in (0,1] 
$,
which approximate 
the initial value functions 
$ 
  \R^d \ni x \mapsto u_d( t, x ) \in \R 
$,
$ d \in \N $, 
without the curse of dimensionality. 
In particular, we observe that the assumption 
in Theorem~\ref{thm:DNN_nonlinear_PDEs} that 
for all 
$ d \in \N $, $ x = (x_1, \dots, x_d) \in \R^d $, 
$ \varepsilon \in (0,1] $, $ t \in [0,T] $ 
it holds that
$ 
  \varepsilon | u_d(t,x) | 
  + \allowbreak 
  | u_d(0,x) \allowbreak - ( \mathcal{R}(\mathfrak{g}_{d,\varepsilon}) )(x) | 
  \le \varepsilon \kappa d^\kappa (1 + \sum_{ i = 1 }^d | x_i |^\kappa ) 
$
assures that 
for all 
$ d \in \N $, $ x = (x_1, \dots, x_d) \in \R^d $, 
$ \varepsilon \in (0,1] $
it holds that
$ 
  | u_d(0,x) \allowbreak - ( \mathcal{R}(\mathfrak{g}_{d,\varepsilon}) )(x) | 
  \le \varepsilon \kappa d^\kappa (1 + \sum_{ i = 1 }^d | x_i |^\kappa ) 
$
and this condition, in turn, ensures 
that for all $ d \in \N $, $ x \in \R^d $ 
it holds that 
$
  ( \mathcal{R}(\mathfrak{g}_{d,\varepsilon}) )(x)
$
converges to 
$
  u_d(0,x)
$
as 
$ \varepsilon $
converges to $ 0 $.

Moreover, we observe that the assumption 
in Theorem~\ref{thm:DNN_nonlinear_PDEs} that 
for all $ d \in \N $, 
$ \varepsilon \in (0,1] $
it holds that 
$
  \mathcal{P}(\mathfrak{g}_{d,\varepsilon}) \le \kappa d^\kappa \varepsilon^{-\kappa} 
$
assures that the number of parameters 
of the neural networks 
$ \mathfrak{g}_{d,\varepsilon} \in {\bf N} $,
$ 
  d \in \N 
$,
$
  \varepsilon \in (0,1] 
$,
grows at most polynomially 
in both the reciprocal $ \varepsilon^{ - 1 } $ 
of the prescribed approximation precision 
$ \varepsilon \in (0,1] $
and the PDE dimension $ d \in \N $.

Furthermore, we note that 
the assumption in Theorem~\ref{thm:DNN_nonlinear_PDEs} 
that for all 
$ d \in \N $, $ x = (x_1, \dots, x_d) \in \R^d $, 
$ \varepsilon \in (0,1] $, $ t \in [0,T] $ 
it holds that 
$ 
  \varepsilon | u_d(t,x) | 
  + \allowbreak 
  | u_d(0,x) \allowbreak - ( \mathcal{R}(\mathfrak{g}_{d,\varepsilon}) )(x) | 
  \le \varepsilon \kappa d^\kappa (1 + \sum_{ i = 1 }^d | x_i |^\kappa ) 
$
demonstrates that 
for all 
$ d \in \N $, $ x = (x_1, \dots, x_d) \in \R^d $, 
$ t \in [0,T] $ 
it holds that 
$ 
  | u_d(t,x) | 
  \le \kappa d^\kappa (1 + \sum_{ i = 1 }^d | x_i |^\kappa ) 
$
and this condition, in turn, ensures 
that the solutions of the PDEs in \eqref{eq:1a} 
grow at most polynomially 
in both the space variable $ x \in \R^d $
and the PDE dimension $ d \in \N $. 
The condition that 
for all 
$ d \in \N $, $ x = (x_1, \dots, x_d) \in \R^d $, 
$ t \in [0,T] $ 
it holds that 
$ 
  | u_d(t,x) | 
  \le \kappa d^\kappa (1 + \sum_{ i = 1 }^d | x_i |^\kappa ) 
$
also ensures that the solutions of the PDEs in \eqref{eq:1a} 
are uniquely described by their initial value 
functions 
$
  \R^d \ni x \mapsto u_d(0,x) \in \R
$,
$ d \in \N $ (cf., e.g., Beck et al.~\cite[Theorem~1.1]{beck2020nonlinear}).

Roughly speaking, \emph{the conclusion} of Theorem~\ref{thm:DNN_nonlinear_PDEs} 
assures that there exist neural networks 
$ \mathfrak{u}_{d,\varepsilon} \in {\bf N} $,
$ d \in \N $, $ \varepsilon \in (0,1] $,
such that for all 
$ d \in \N $, $ \varepsilon \in (0,1] $ 
it holds that the $ L^2 $-approximation error 
$
  [ 
    \int_{ [0,1]^d } 
      | 
        u_d(T,x) - 
        ( \mathcal{R}( \mathfrak{u}_{ d, \varepsilon } ) )(x) 
      |^2 
    \,
    dx 
  ]^{ \nicefrac{ 1 }{ 2 } } 
$
between the exact solution 
$
  u_d(T,x)
$
of the PDE 
and its neural network approximation 
$
  ( \mathcal{R}( \mathfrak{u}_{ d, \varepsilon } ) )(x) 
$
is smaller or equal than 
the prescribed approximation accuracy 
$ \varepsilon $ 
while the numbers  
$
  \mathcal{P}( \mathfrak{u}_{d,\varepsilon} ) 
$,
$ d \in \N $,
$ \varepsilon \in (0,1] $,
of parameters 
of the approximating neural networks 
$ \mathfrak{u}_{d,\varepsilon} \in {\bf N} $,
$ d \in \N $, $ \varepsilon \in (0,1] $, 
grow at most polynomially 
in both 
the PDE dimension $ d \in \N $ 
and the reciprocal $ \varepsilon^{ - 1 } $ 
of the prescribed approximation accuracy $ \varepsilon $. 
We note that Theorem~\ref{thm:DNN_nonlinear_PDEs} 
above is a neural network approximation result 
for the solutions 
of the PDEs in \eqref{eq:1a} at the final time $ T $ 
on the $ d $-dimensional hypercube $ [0,1]^d $ 
but the more general neural network approximation result 
in Hutzenthaler et al.~\cite[Theorem~4.1]{hutzenthaler2020proof}
also provides neural networks approximations 
for solutions 
of PDEs 
on more general space regions.

In the next step we add some words on 
\emph{the strategy 
of the proofs} of 
Theorem~\ref{thm:DNN_nonlinear_PDEs} above 
and 
Theorem~1.1 in Hutzenthaler et al.~\cite{hutzenthaler2020proof}, 
respectively. Even though 
Theorem~\ref{thm:DNN_nonlinear_PDEs} above
and 
Theorem~1.1 in Hutzenthaler et al.~\cite{hutzenthaler2020proof}, respectively, 
are purely deterministic neural network 
approximation results, 
the proofs of 
Theorem~\ref{thm:DNN_nonlinear_PDEs} above 
and 
Theorem~1.1 in Hutzenthaler et al.~\cite{hutzenthaler2020proof}, 
respectively, are strongly based on 
probabilistic arguments on a suitable 
artificial probability space. 
In particular, the proofs of 
Theorem~\ref{thm:DNN_nonlinear_PDEs} above 
and 
Theorem~1.1 in Hutzenthaler et al.~\cite{hutzenthaler2020proof}, 
respectively, 
employ the fact in the following elementary lemma. 
\begin{lemma}
\label{lem:elementary}
Let $ \varepsilon \in (0,\infty) $, 
let $ ( \Omega, \mathcal{F}, \P ) $ be a probability space, 
and 
let $ E \colon \Omega \to \R $ be a random variable 
with $ ( \E[ | E |^2 ] )^{ 1 / 2  } \leq \varepsilon $. 
Then there exists $ \omega \in \Omega $ such that 
$ | E( \omega ) | \leq \varepsilon $.
\end{lemma}
The elementary statement in Lemma~\ref{lem:elementary} follows, e.g., 
from Grohs et al.~\cite[Proposition~3.3]{grohs2018proof}.
Lemma~\ref{lem:elementary} is employed in the proofs of 
Theorem~\ref{thm:DNN_nonlinear_PDEs} above 
and 
Theorem~1.1 in Hutzenthaler et al.~\cite{hutzenthaler2020proof}, 
respectively, 
to construct an appropriate random realization 
with desired approximation properties 
on a suitable artifical probability space. 
To make it more concrete, 
the proofs of 
Theorem~\ref{thm:DNN_nonlinear_PDEs} above 
and 
Theorem~1.1 in Hutzenthaler et al.~\cite{hutzenthaler2020proof}, respectively,  
consist, roughly speaking, of the following 
four steps (cf., e.g., \cite[Section~1]{jentzen2018proofarxiv1809}):
\begin{enumerate}[(I)]
\item 
\label{item:I}
First, appropriate random neural networks are constructed 
on a suitable artificial probability space. 
These appropriate neural networks 
are random in the sense that the weights and the biases of these 
neural networks are random variables instead of deterministic real numbers. 
The random neural networks are appropriately constructed with the aim 
to appropriately approximate the solutions of the PDEs in \eqref{eq:1a}.
\item 
\label{item:II}
Second, it is proved that the realization functions of 
these random neural networks are in a suitable root mean square sense 
close to the solutions of the PDEs in \eqref{eq:1a} at the final time $ T $. 
\item 
\label{item:III}
Third, it is proved that the numbers of parameters 
of these random neural networks 
grow at most polynomially in both 
the reciprocal $ \varepsilon^{ - 1 } $ of the 
prescribed approximation accuracy $ \varepsilon \in (0,1] $ 
and the PDE dimension $ d \in \N $. 
Here the approximation accuracy is measured in a suitable root mean square sense 
according to item~\eqref{item:II} above.
\item 
Fourth, Lemma~\ref{lem:elementary} is applied 
to suitable error random variables, which describe 
certain $ L^2 $-errors between the realization functions 
of the constructed random neural networks 
(see item~\eqref{item:II} above)
and the exact solutions of the PDEs in \eqref{eq:2a} at the final time $ T $
(cf.\ \eqref{eq:2a} in Theorem~\ref{thm:DNN_nonlinear_PDEs} above), 
to obtain the existence of a realization on 
the artifical probability space such that 
the error random variables evaluated 
at this realization are smaller or equal than 
the prescribed approximation accuracy $ \varepsilon $. 
Combining the existence of such a realization 
on the artificial probability space 
with item~\eqref{item:III} above 
then completes 
the proofs of 
Theorem~\ref{thm:DNN_nonlinear_PDEs} above 
and 
Theorem~1.1 in 
Hutzenthaler et al.~\cite{hutzenthaler2020proof}, 
respectively. 
\end{enumerate}
Let us also add a few comments 
on the way how the appropriate random neural networks 
in item~\eqref{item:I} above are designed 
and on the way how the statements sketched in 
items~\eqref{item:II}--\eqref{item:III} above are proved. 
The main tool for items~\eqref{item:I}--\eqref{item:III} above 
are MLP approximations 
(cf.\ Section~\ref{sec:MLP} above). 
More formally, 
the random neural networks in item~\eqref{item:I} above 
are designed so that their realization functions 
coincide with suitable MLP approximations 
and the statement in item~\eqref{item:II} above 
is then proved by employing suitable root mean square error 
estimates for MLP approximations 
(cf.\ Hutzenthaler et al.~\cite[Theorem~3.5]{Hutzenthaleretal2018arXiv} and Theorem~\ref{thm:MLP_analysis} above)
and the statement in item~\eqref{item:III} above 
is then proved by employing suitable cost estimates 
for neural networks and MLP approximations 
(cf.\ Hutzenthaler et al.~\cite[Sections~3.2--3.3]{hutzenthaler2020proof}).

\section{Conclusion}
The progress reviewed here  has opened up a host of  new possibilities, both in theory and applications.
In applications, it has been proved effective in finance, such as the pricing of financial derivatives~\cite{wang2018deep,BeckerCheridito2019,BeckerCheriditoJentzen2019,becker2020pricing} and credit valuation adjustment~\cite{gnoatto2020deep}. 
It also opens up new possibilities in control theory, an area that has long been hindered by the curse of dimensionality problem.
In fact, it is likely that control theory will be among the areas most impacted by the kind of ideas reviewed here.

Another interesting new problem is the mathematical study of high dimensional PDEs.  The fact that we can compute their solutions rather efficiently even in very high dimensions means that the complexity of these solutions should not be very high.  Can we quantify this in some way?
In low dimensions, a  huge amount of effort has gone into studying the regularity of solutions of PDEs.  It seems that regularity is not the most important issue in high dimensions.  Rather, it is the complexity that is more relevant. It would be very interesting to develop a complexity-based PDE theory in high dimensions.  It is worth mentioning that in low dimensions, regularity is also
a measure of complexity:  The efficiency of approximating a target function by certain approximation scheme, say piecewise polynomial approximation, is often measured by the regularity of the target function.

An interesting topic untouched in this review is reinforcement learning.  Formulated with the language we use here, reinforcement learning is  all about solving the Bellman equation for the underlying Markov decision process \cite{sutton2018reinforcement}.
However, in contrast to the ideas reviewed here, which make heavy use of the underlying model, reinforcement learning makes minimum use of the specifics of the model. At this moment, it is still quite unclear what the relative merits are between the ideas reviewed here and those of reinforcement learning.  This is undoubtedly an interesting area for further work.

\section*{Acknowledgement}
The third author acknowledges funding by the Deutsche Forschungsgemeinschaft (DFG, German Research Foundation) under Germany’s Excellence Strategy EXC 2044-390685587, Mathematics M\"{u}nster: Dynamics-Geometry-Structure.

\bibliographystyle{acm}
\bibliography{ref}

\end{document}